\documentclass[11pt, ,bezier,amstex]{article}
\usepackage{amssymb,amsmath,latexsym}
\usepackage{epsfig}

\usepackage{epic,eepic,wrapfig,color}

\newtheorem{thm}{Theorem}[section]
\newtheorem{cor}[thm]{Corollary}
\newtheorem{lemma}[thm]{Lemma}
\newtheorem{prop}[thm]{Proposition}
\newtheorem{defn}[thm]{Definition}

\newtheorem{remark}[thm]{Remark}
\newtheorem{example}[thm]{Example}
\numberwithin{equation}{section}

\def\pf{{\medskip\noindent {\bf Proof. }}}
\def\qed{{\hfill $\Box$ \bigskip}}
\def\R{{\mathbb R}}
\def\P{{\mathbb P}}
\def\E{{\mathbb E}}
\def\1{{\bf 1}}

\def\sJ {{\cal J}}  \def\sL {{\cal L}}

\def\bP {{\mathbb P}}  \def\bR {{\mathbb R}}

\def\R {{\mathbb R}}

\def\nn{\nonumber}

\def\wt{\widetilde}
\def\wh{\widehat}

\def\E{{\mathbb E}}
\def\P{{\mathbb P}}

\def\bea{\begin{align*}}
\def\eea{\end{align*}}
\def\bee{\begin{equation}}
\def\eee{\end{equation}}

\def\eps{\varepsilon}
\def\vp{\varphi}

\def\wh{\widehat}

\makeatletter
\@addtoreset{equation}{section}

\makeatother

\begin{document}
\bibliographystyle{plain}

\title{\Large \bf
Dirichlet Heat Kernel Estimates for  \\ Rotationally Symmetric L\'evy processes}

\author{{\bf Zhen-Qing Chen}\thanks{Research partially supported
by NSF Grant DMS-1206276 and NNSFC Grant 11128101.}, \quad {\bf Panki Kim}\thanks{This work was supported by Basic Science Research Program through the National Research Foundation of Korea(NRF)
grant funded by the Korea government(MEST)
(2012-0000940). } \quad and  \quad {\bf Renming Song}\thanks{Research supported in part by a grant from the Simons
Foundation (208236).}}
\date{(March 26, 2013)}

\maketitle

\begin{abstract}
In this paper, we consider a large class of
purely discontinuous rotationally symmetric L\'evy
processes.
We establish sharp two-sided estimates for
the transition densities of such processes
killed upon leaving an open set $D$.
When $D$ is a $\kappa$-fat open set, the sharp two-sided estimates
are given in terms of surviving probabilities
and the global transition density of the L\'evy process.
When $D$ is a $C^{1, 1}$ open set and
the L\'evy exponent of the process is given by $\Psi(\xi)= \phi(|\xi|^2)$
with $\phi$ being a complete Bernstein function
satisfying a mild growth condition at infinity,
our two-sided estimates are explicit in terms of $\Psi$, the distance
function to the boundary of $D$ and the jumping kernel of
$X$, which give an affirmative answer to the conjecture posted in
\cite{CKS6}.
Our results are the first sharp two-sided
Dirichlet heat kernel estimates for
a large class of symmetric L\'evy processes
with general L\'evy exponents.
We also derive an explicit lower bound estimate for symmetric L\'evy processes
on $\R^d$ in terms of their L\'evy exponents.
\end{abstract}


\bigskip

\section{Introduction}\label{s:I}
Due to their importance in theory and applications, fine potential
theoretical properties of L\'evy processes have been under
intense study recently. The transition density $p(t, x, y)$ of a
L\'evy process is the heat kernel of the generator of the process.
However, the transition density  (if it exists) of a general L\'evy process
rarely admits an explicit expression.
Thus obtaining sharp estimates on $p(t, x, y)$ is a fundamental problem both
in probability theory and in analysis.

The generator of a discontinuous L\'evy process is an integro-differential
operator and so it is a non-local operator. Recently, quite a few people in PDE
are interested in problems related to non-local operators; see, for example,
\cite{CS1, CS2, CSS, S} and the references therein.

   When $X$ is a symmetric diffusion on $\R^d$ whose infinitesimal generator
  is a uniformly elliptic and bounded divergence form operator,
  it is well-known that
  $p(t, x, y)$ enjoys the celebrated Aronson's Gaussian type estimates.
 When $X$ is a pure jump symmetric process on $\R^d$,
sharp estimates on $p(t, x, y)$ have been
studied in \cite{CKK, CKK2, CKK3, CK, CK2} recently,
which can be viewed as the counterpart of Aronson's estimates
for non-local operators.

Due to the complication near the boundary, two-sided estimates for the
transition densities of
discontinuous L\'evy processes
killed upon leaving  an open set $D$
(equivalently, the Dirichlet heat kernels) have been established
very recently for a few particular processes only.
The first of such estimates is obtained in \cite{CKS}, where
we succeeded in establishing sharp two-sided estimates for the heat
kernel of the fractional Laplacian
$\Delta^{\alpha /2}:=-(-\Delta)^{\alpha /2}$ with zero
exterior condition on $D^c$ (or equivalently, the transition density
of the killed symmetric  $\alpha$-stable process) in any $C^{1, 1}$
open set $D$. The approach developed in \cite{CKS} provides a road map
for establishing sharp two-sided heat kernel estimates of
other jump processes in open subsets of $\R^d$. The ideas of \cite{CKS}
were adapted to establish sharp
two-sided heat kernel estimates of relativistic stable processes and
mixed stable processes in
$C^{1, 1}$ open subsets of $\R^d$ in \cite{CKS2, CKS4} respectively.
In all these cases, the characteristic exponents of these L\'evy processes
admit explicit expressions, the boundary decay rates of the Dirichlet heat kernels are
suitable powers of the distance to the boundary.
On the other hand, a Varopoulos type two-sided
Dirichlet heat kernel estimate  of symmetric stable
processes in $\kappa$-fat open sets
was derived in \cite{BGR}; this type of estimates is
expressed in terms of surviving probabilities
and the global transition density of
the symmetric stable process.

The objective of this paper is to establish sharp two-sided estimates on
the transition density $p_D(t, x, y)$ for a large class of
purely discontinuous rotationally symmetric  L\'evy processes. Unlike
the cases considered in
\cite{BGR, CKS,  CKS4, CKS2},
the characteristic exponents of the symmetric L\'evy processes considered
in this paper are quite general, satisfying only certain mild growth condition
at $\infty$.
Moreover, the boundary decay rate of $p_D(t, x, y)$ is no longer
some power of the distance to the boundary.
The analysis of the precise boundary
behavior of $p_D(t, x, y)$ is
 quite challenging and delicate. The main tools to obtain the precise boundary
behavior of $p_D(t, x, y)$ are two versions of the boundary Harnack principle obtained in \cite{KM, KSV5}.
In this paper we combine the approaches developed in \cite{BGR, CKS} with these boundary Harnack
principles
 to obtain sharp two-sided estimates for
$p_D(t, x, y)$, which cover the main results in
\cite{BGR, CKS, CKS4, CKS2} and much more.

Suppose that $S=(S_t: t\ge 0)$ is a subordinator with Laplace exponent $\phi$, that is, $S$ is
a nonnegative L\'evy process with $S_0=0$ and
$ \E\left[e^{-\lambda S_t}\right]=e^{-t\phi(\lambda)}$ for every $t, \lambda >0$.
The function $\phi$ can be written in the form
\begin{equation}\label{e:LK}
\phi(\lambda)=b\lambda + \int^\infty_0(1-e^{-\lambda t})\mu(dt),
\end{equation}
where $b\ge 0$ and $\mu$ is a measure on $(0, \infty)$ satisfying $\int^\infty_0
(1\wedge t)\mu(dt)<\infty$.
The constant $b$
is called the drift of
the subordinator and $\mu$
the L\'evy measure of the subordinator (or of $\phi$).
The function $\phi$ is a Bernstein function, i.e., it is
$C^{\infty}$, positive and $(-1)^{n-1} D^n \phi\ge 0$ for all $n\ge 1$.
In particular, since $\phi (0)=0$ and $\phi ''\leq 0$,
the Bernstein function $\phi$ has the property that
\begin{equation}\label{e:Berall}
\phi(\lambda r )\le \lambda\phi(r)
\qquad \text{ for all }
\lambda \ge 1 \hbox{ and }  r >0.
\end{equation}

The Laplace exponent  $\phi$ is said to be
a complete Bernstein function if  the L\'evy measure
$\mu$ of $\phi$ has a
completely monotone density $\mu(t)$, i.e., $(-1)^n D^n\mu\ge 0$ for every non-negative integer $n$.
For basic results on complete
Bernstein functions,
we refer the reader to \cite{SSV}.

Throughout this paper,  we assume that
$\phi$ is a complete Bernstein function  satisfying the following
growth condition at infinity (see \cite{Z}):

\medskip
\noindent
{\bf (A):}
There exist constants $ \delta_1, \delta_2 \in (0,1)$, $a_1\in (0, 1)$, $a_2\in (1, \infty)$
and $R_0>0$ such that
\begin{eqnarray*}
a_1\lambda^{\delta_1} \phi(r) \leq
\phi(\lambda r) \le a_2 \lambda^{\delta_2} \phi(r) \quad
&\hbox{for } \lambda \ge 1 \hbox{ and }  r \ge R_0.
\end{eqnarray*}
Note that it follows from
the upper bound condition in {\bf (A)}
that $\phi$ has no drift.

Let $W=(W_t:\, t\ge 0)$ be a Brownian motion in $\R^d$ independent of the subordinator $S$.
The subordinate Brownian motion $Y=(Y_t:\, t\ge 0)$ is defined by $Y_t:=
W_{S_t}$,
which is a rotationally symmetric L\'evy process with L\'evy exponent $\phi
( |\xi|^2)$. The infinitesimal generator of $Y$ is $\sL^Y := - \phi (-\Delta)$.
 Here and below for a function $\psi$ on $[0, \infty)$, $\psi (-\Delta)$ is defined as a pseudo differential operator in terms of Fourier transform;
that is, $\wh {\psi (-\Delta) f} (\xi) := - \psi (|\xi|^2) \wh f(\xi)$, where 
$\wh f$ is the Fourier transform of a function $f$ on $\R^d$.
 It is known that
the L\'evy measure of the process $Y$ has a density given by
$J(x)=j(|x|)$ where
\begin{equation}\label{e:ld4s}
j(r):=\int^{\infty}_0(4\pi t)^{-d/2}e^{-r^2/(4t)}\mu(t)dt, \qquad r>0.
\end{equation}
 Note that the function
$r\mapsto j(r)$ is continuous and decreasing
on $(0, \infty)$.

We will assume that $X$ is a purely discontinuous rotationally
symmetric L\'evy process with L\'evy exponent
$\Psi(\xi)$.  Because of rotational symmetry, the function $\Psi$ depends on
$|\xi|$ only, and by a slight abuse of notation we write $\Psi(\xi)=\Psi(|\xi|)$.
The infinitesimal generator of $X$ is $\sL^X:= - \Psi (\sqrt{-\Delta})$.
We further assume that the L\'evy measure of $X$ has a
density with respect to the Lebesgue measure on $\R^d$,
which is denoted by
$J_X(x, y)=J_X(x-y)=j_X(|y-x|)$.
That is,
$$
\E_x\left[e^{i\xi\cdot(X_t-X_0)}\right]=e^{-t\Psi(|\xi|)}
\quad \quad \mbox{ for every } x\in \R^d \mbox{ and } \xi\in \R^d,
$$
with
\begin{equation}\label{e:psi}
\Psi(|\xi|)= \int_{\R^d}(1-\cos(\xi\cdot x))J_X(x)dx.
\end{equation}
We assume that $j_X (r)$ is continuous on $(0, \infty)$
and that there is a constant $\gamma>1$ such that
\begin{equation}\label{e:psi1}
\gamma^{-1} j(r )\le J_X(r) \le \gamma j(r)
 \quad \mbox{for all } r>0 .
\end{equation}
This implies that $J$ and $J_X$ are comparable.
Clearly \eqref{e:psi1} also implies that
\begin{equation}\label{e:psi2}
\gamma^{-1} \phi(|\xi|^2)\le \Psi(|\xi|) \le \gamma \phi(|\xi|^2)
 \quad \mbox{for all } \xi\in \R^d\, .
\end{equation}

We remark that under the above assumptions, $X$ does not need to be a
subordinate Brownian motion  because $j_X$ does not need to be monotone.
For example, choose $\eps>0$ such that
$2^{-1}j(1)< j(1+\eps)$ and a continuous function $h$
with $h(1)=\frac12 j(r)$, $h(1+\eps)=0$ and $0\le h(r) \le \frac12 j(r)$
for all $r>0$.
Then $j_X(r):=j(r)-h(r)$ is not monotone and its corresponding
L\'evy process $X$ (through L\'evy exponent \eqref{e:psi})
is not a subordinate Brownian motion.

Under the above setup, $X$ has a continuous transition density
$p(t, x, y)$ with respect to the Lebesgue measure on $\R^d$ (see \cite{CKK2}).
Clearly, $p(t, x, y)$ is a function depending only on $t$ and $|x-y|$, and so,
by an abuse of notation, we also denote $p(t, x, y)$ by $p(t, |x-y|)$.
For every open subset $D\subset \bR^d$, we denote by  $X^D$  the
subprocess of $X$ killed upon exiting $D$. It is known (see \cite{CKK2}) that $X^D$ has a
transition density $p_D(t, x, y)$, with respect to the Lebesgue
measure, which is jointly locally H\"older continuous.
Note that $p_D (t, x, y)$ is the fundamental solution for
$\sL^X=-\Psi (\sqrt{-\Delta})$ in $D$ with zero exterior condition
 and so it can also
be called the Dirichlet heat kernel of $\sL^X$ in $D$.
The purpose of this paper is to establish
sharp two-sided estimates on $p_D(t, x, y)$.
The following two conditions will be needed for some of the results
in this paper when $D$ is unbounded.

\medskip
\noindent {\bf (B)}:
There exist constants $C_1>0$ and $C_2 \in (0, 1]$ such that
$$ p(t, u) \leq C_1 p(t, C_2 r )
\qquad \hbox{for } t\in (0, 1] \hbox{ and }
 u\geq r>0.
$$

\medskip
\noindent {\bf (C)}:
There exist constants $C_3>0$ and $C_4 \in (0, 1]$ such that
$$ p(t, r) \leq
C_3 t j (C_4 r)
\qquad \hbox{for } t\in (0, 1] \hbox{ and }
r>0.
$$

\bigskip
Throughout this paper we will use $\Phi$ to denote the function
\bee \label{e:1.9}
\Phi(r)=\frac1{\phi(r^{-2})}, \qquad r>0.
\eee
Note that in particular it follows from \eqref{e:Berall} that
\bee\label{e:1.8}
\Phi (2r)=\frac{1}{\phi (r^{-2}/4)}\leq \frac{1}{\phi^{-1}(r^{-2})/4}
=4\Phi (r) \qquad \hbox{for every } r>0.
\eee
The inverse function of $\Phi$ will
be denoted by the usual notation $\Phi^{-1}(r)$.
Here and in the following, for $a, b\in \bR$,  $a\wedge b:=\min \{a, b\}$
and $a\vee b:=\max\{a, b\}$.

\bigskip

\begin{remark}\label{R:1.1} \rm
\begin{description}  \item{(i)}
The condition {\bf (B)} is pretty mild. When $X$ is a rotationally
symmetric L\'evy process such that $r\mapsto j_X(r)$ is decreasing,
condition {\bf (B)} holds for all $t>0$ (see \cite[Proposition]{W}).
In particular it holds for
 all subordinate Brownian motions with
$C_1=C_2=1$.
In this special case, we can also see this using the following elementary argument:
When $X$ is a subordinate Brownian motion,
$p(t, r)=\E [ p_0(S_t, r)]$, where
$p_0(t, |x-y|)$
is the transition
density of the Brownian motion $W$.
It follows immediately that $p(t, r)$
is   decreasing in $r$ and so {\bf (B)} holds
with $C_1=C_2=1$.

\item{(ii)} Under condition {\bf (A)},
condition {\bf (B)} is weaker than {\bf (C)}.
Under condition {\bf (A)}, we will show in this paper that
there exists $c>1$ such that
$$
c^{-1} \big(\Phi^{-1} (t)^d \wedge t j(r) \big)
\le p(t, r)\leq c (\Phi^{-1} (t))^d
\qquad \hbox{for } t\in (0, 1] \hbox{ and } r>0
$$
(see Proposition \ref{p:new1} and Theorem \ref{t:globalhke} below).
Thus condition {\bf (C)} amounts to say
that  there exist constants $c\geq 1$ and $C_4\in (0, 1]$ such that
for $t\in (0, 1]$ and $r>0$,
\begin{equation}\label{e:1.7}
 c^{-1}  \big( \Phi^{-1} (t)^d \wedge t j(r) \big)
 \leq p(t, r)   \leq c
 \big( \Phi^{-1} (t)^d \wedge t j(C_4 r) \big) .
\end{equation}
Since $j(r)$ is a decreasing function in $r$, \eqref{e:1.7}
implies that condition {\bf (B)} holds with $C_1=c^2$ and $C_2=C_4$.

\item{(iii)} Assume that condition {\bf (A)} holds.
It follows from
\cite{CKK2, CKK3, CK2} that, for every $R>0$, there is a constant
 $c=c(T, R, \gamma, \phi)>1$ so that \eqref{e:1.7} holds
for $(t, r)\in (0, 1]\times (0, R]$ with $C_4=1$. (See Proposition \ref{p:new1} below.)
So the assertions in conditions {\bf (B)} and {\bf (C)}
are always satisfied for $0<r \le u\leq R$.

\item{(iv)} By
\cite{CKK3, CK2},
 $$
 c^{-1} \left( ( \Phi^{-1}(t))^{-d} \wedge t j  (C_4^{-1} r)\right)
\leq p(t, r) \leq c \left( ( \Phi^{-1}(t))^{-d} \wedge t j  (C_4 r)\right)
\quad  \hbox{for } (t, r)\in (0, 1]\times (0, \infty),
$$
and consequently conditions {\bf (B)} and {\bf (C)},
 hold  for a large class of discontinuous processes
including mixed stable-like processes
(with $C_4=1$) and relativistic stable-like processes (with $C_4=1$,
see \cite[Theorem 4.1]{CKS2}).
\qed
\end{description}
\end{remark}

Before stating the main results of this paper,
we need first to set up some notations.
Let $d \ge 1 $. We denote the Euclidean distance between
$x$ and $y$ in $\R^d$ by  $|x-y|$ and
denote by $B(x, r)$ the open
ball centered at $x\in \bR^d$ with radius $r>0$;
for any two positive functions $f$ and $g$,
$f\asymp g$ means that there is a positive constant $c\geq 1$
so that $c^{-1}\, g \leq f \leq c\, g$ on their common domain of
definition;
for any open $D\subset\bR^d$ and $x\in D$,
${\rm diam}(D)$ stands for the diameter of $D$ and
$\delta_D(x)$ stands for the Euclidean distance between
$x$ and $D^c$.

\begin{defn}\label{def:UB}
  Let $0<\kappa\leq 1$.
  We say that a open set $D$ is $\kappa$-fat if
  there is $R_1>0$ such that for all $x\in \overline{D}$ and all $r\in
  (0,R_1]$, there is a ball
  $B(A_r(x), \kappa r)
  \subset D\cap B(x,r)$. The pair $(R_1, \kappa)$ is called the characteristics of
the $\kappa$-fat open set $D$.
\end{defn}

The following factorization of the Dirichlet heat kernel is the first main result of this paper.
Recall that $C_1$ and $C_2$ are the constants in condition {\bf (B)}.

\begin{thm}\label{thm:oppz}
Let $X$ be a purely discontinuous rotationally symmetric L\'evy
process with L\'evy exponent $\Psi$ and L\'evy density
$J_X$ satisfying \eqref{e:psi2} and \eqref{e:psi1} respectively,
 where the complete Bernstein function $\phi$ satisfies {\bf (A)}.
  Suppose that  $D$ is a $\kappa$-fat open set with
characteristics $(R_1, \kappa)$.
    \begin{description}
\item{\rm (i)}
For every $T>0$,
there exists  $c_1=c_1( R_1, \kappa,
\gamma,  T, d,  \phi)>0$
 such that for $0<t\leq T$, $x,y\in D$,
   \begin{equation}   \label{eq:get2}
 p_D(t, x, y)\geq  c_1  \P_{x}(\tau_{D}>t)\P_{y}(\tau_{D}>t)
 \left( \Phi^{-1}(t)^{-d} \wedge t J(x, y)\right).
  \end{equation}
\item{\rm (ii)} If $D$ is unbounded, we assume in addition that condition {\bf (B)} holds.
For every $T>0$,
there exists  $c_2=c_2(C_1, C_2, R_1, \kappa, T, d, \gamma, \phi)>0$
 such that for $0<t\leq T$, $x,y\in D$,
\begin{equation} \label{eq:get}
     {p_{D}(t, x, y)}
    \le
    c_2 \P_{x}(\tau_{D}>t)
    \P_{y}(\tau_{D}>t){p(t,C_5 x,C_5 y)},
  \end{equation}
 where $C_5 = C_2^2/4$.

  \item{\rm (iii)} Suppose in addition that $D$ is bounded.
Then there exists
$c_3=c_3(\text{diam}(D), R_1, \kappa,  d, \gamma, \phi)>1$
so that for all $(t, x, y)\in [3, \infty)\times
D\times D$,
$$
c_3^{-1}\, \P_x(\tau_D > 1)\,\P_y(\tau_D > 1)\, e^{-t\lambda_1}
\,\leq\,  p_D(t, x, y) \,\leq\, c_3\,  \P_x(\tau_D > 1)\,\P_y(\tau_D > 1) \,e^{-t\lambda_1},
$$
where $-\lambda_1<0$ is the largest eigenvalue of
the generator of $X^D$.
\end{description}\end{thm}

When $X$ is a rotationally symmetric $\alpha$-stable process in $\R^d$,
that is, when $\Psi (\xi)=|\xi|^\alpha$ for some $\alpha \in (0, 2)$,
parts (i) and (ii) of Theorem \ref{thm:oppz} are proved in \cite{BGR}.

Recall that $C_3$ and $C_4$ are the constants in condition {\bf (C)}.
Combining Theorem \ref{thm:oppz}(i)--(ii) and Remark \ref{R:1.1}(ii),
we have the following corollary.

\begin{cor}\label{c:main1}
Let $X$ be a purely discontinuous rotationally symmetric L\'evy
process with L\'evy exponent $\Psi$ and L\'evy density
$J_X$ satisfying \eqref{e:psi2} and \eqref{e:psi1} respectively,
 where the complete Bernstein function $\phi$ satisfies {\bf (A)}.
  Suppose that  $D$ is a $\kappa$-fat open set with
characteristics $(R_1, \kappa)$.
 If $D$ is unbounded, we assume in addition that condition {\bf (C)} holds.
For every $T>0$,
there exist  $c_1=c_1(R_1, \kappa, T, d, \gamma, \phi)>0$
and $c_2=c_2(C_3, C_4, R_1, \kappa, T, d, \gamma, \phi)>0$
 such that for $0<t\leq T$, $x,y\in D$,
 \begin{align*}
   &  c_1\P_{x}(\tau_{D}>t)\P_{y}(\tau_{D}>t)
 \left( \Phi^{-1}(t)^{-d} \wedge t j(|x- y|)\right)\\
 \le \, & {p_{D}(t, x, y)}
    \le c_2
 \P_{x}(\tau_{D}>t) \P_{y}(\tau_{D}>t)
 \left( \Phi^{-1}(t)^{-d} \wedge t j(C_6|x- y|)\right),
\end{align*}
 where $C_6 = C_4^3/4$.
\end{cor}

The second main result of this paper is on  explicit sharp
Dirichlet heat kernel estimates
for subordinate Brownian motions in $C^{1,1}$ open sets.
So in the remainder of this section, we assume that
$X=Y$, a subordinate Brownian motion with L\'evy exponent $\Psi(\xi)=
\phi(|\xi|^2)$.

Recall that an open
set $D$ in $\bR^d$ (when $d\ge 2$) is said to be a (uniform)
$C^{1,1}$ open set if there exist a localization radius $R_2>0$ and
a constant $\Lambda>0$ such that for every $z\in\partial D$, there
exist a $C^{1,1}$-function $\psi=\psi_z: \bR^{d-1}\to \bR$
satisfying $\psi (0)= 0$, $\nabla
\psi (0)=(0, \dots, 0)$, $\| \nabla
\psi  \|_\infty
\leq \Lambda$, $| \nabla
\psi (x)-\nabla
\psi (z)| \leq \Lambda
|x-z|$, and an orthonormal coordinate system $CS_z$ with its origin
at $z$ such that
$$
B(z, R_2)\cap D=\{ y= (\wt y, \, y_d) \mbox{ in } CS_z: |y|< R_2,
y_d > \psi (\wt y) \}.
$$
The pair $(R_2, \Lambda)$ is called the characteristics of the
$C^{1,1}$ open set $D$. Note that a $C^{1,1}$ open set $D$ with
characteristics $(R_2, \Lambda)$ can be unbounded and
disconnected; the distance between two distinct components of $D$ is
at least $R_2$.
By a $C^{1,1}$ open set in $\bR$ we mean an open set which can be
written as the union of disjoint intervals so that the minimum of
the lengths of all these intervals is positive and the minimum of
the distances between these intervals is positive.

Here is the second main result of this paper,
 which gives
an affirmative answer to the Conjecture posed in \cite{CKS6}.
In view of Remark \ref{R:1.1}, it extends the main results of
\cite{CKS, CKS4, CKS2}.
Recall  that
condition {\bf (B)} holds with $C_1=C_2=1$ for any
subordinate Brownian motion.

\begin{thm}\label{t:main}
Suppose that $X$ is a
subordinate Brownian motion with L\'evy exponent $\Psi(\xi)=
\phi(|\xi|^2)$ with $\phi$ being a complete Bernstein function satisfying
condition {\bf (A)}.
Let $D$ be a $C^{1,1}$ open subset of $\bR^d$
with characteristics
$(R_2, \Lambda)$.
\begin{description}
\item{\rm (i)}
For every $T>0$, there exists
$c_1=c_1(R_2, \Lambda, T, d, \phi)>0$
such that for all $(t, x, y) \in (0, T]\times D\times D$,
\begin{eqnarray*}
p_D (t, x, y)\geq  c_1 \left(1\wedge \frac{\Phi (\delta_D(x))}t \right)^{1/2}
\left(1\wedge \frac{\Phi (\delta_D(y))}t \right)^{1/2}
\left(\Phi^{-1}(t)^{-d}\wedge tJ(x, y)\right).
\end{eqnarray*}

\item{\rm (ii)}
For every $T>0$, there exists
$c_2=c_2(
 R_2, \Lambda, T, d, \phi)>0$
such that for all $(t, x, y) \in (0, T]\times D\times D$,
$$
  p_D(t, x, y)
 \leq
 c_2\left(1\wedge \frac{\Phi (\delta_D(x))}t \right)^{1/2} \left(1\wedge
 \frac{\Phi (\delta_D(y))}t \right)^{1/2}
p(t, |x-y|/4).$$

\item{\rm (iii)} Suppose in addition that $D$ is bounded.
For every $T>0$, there exists
$c_3\ge 1$ depending only on
$\text{diam}(D),
 \lambda,
R_2, \Lambda, d, \phi$ and $T$
so that for all $(t, x, y)\in [T, \infty)\times
D\times D$,
$$
c_3^{-1}\, e^{-\lambda_1 t}
\sqrt{\Phi(\delta_D (x))} \sqrt{\Phi(\delta_D (y))}
\,\leq\,  p_D(t, x, y) \,\leq\, c_3\, e^{-\lambda_1 t}\,
\sqrt{\Phi(\delta_D (x))} \sqrt{\Phi(\delta_D (y))},
$$
where $-\lambda_1<0$ is the largest eigenvalue of
the generator of $X^D$.
\end{description}
\end{thm}

Recall  that $C_3$ and $C_4$ are the constants in condition {\bf (C)}.
Combining Theorem \ref{t:main}(i)--(ii) and Remark \ref{R:1.1}(ii), we have the following corollary.

\begin{cor}\label
{c:main2}
Suppose that $X$ is a
 subordinate Brownian motion
 with L\'evy exponent $\Psi(\xi)=
\phi(|\xi|^2)$ with $\phi$ being a complete Bernstein function satisfying
condition {\bf (A)}.
Let $D$ be a $C^{1,1}$ open subset of $\bR^d$
with characteristics
$(R_2, \Lambda)$.
 If $D$ is unbounded, we assume in addition that condition {\bf (C)} holds.
For every $T>0$,
there exist
$c_1=c_1(R_2, \Lambda, T, d, \phi)>0$ and
$c_2=c_2(C_3, C_4, R_2, \Lambda, T, d, \phi)>0$
 such that for $0<t\leq T$, $x,y\in D$,
 \begin{align*}
     &  c_1\left(1\wedge \frac{\Phi (\delta_D(x))}t \right)^{1/2} \left(1\wedge
 \frac{\Phi (\delta_D(y))}t \right)^{1/2}
 \left( \Phi^{-1}(t)^{-d} \wedge t j(|x- y|)\right)\\
 \le \, & {p_{D}(t, x, y)}
    \le c_2 \left(1\wedge \frac{\Phi (\delta_D(x))}t \right)^{1/2} \left(1\wedge
 \frac{\Phi (\delta_D(y))}t \right)^{1/2}\left( \Phi^{-1}(t)^{-d} \wedge t j(C_4|x- y|/4)\right).
\end{align*}
\end{cor}

When $X$ is a
rotationally  symmetric $\alpha$-stable process in $\R^d$,
Theorem  \ref{t:main} is first established in \cite{CKS}.
Sharp two-sided Dirichlet heat kernel estimates in $C^{1,1}$ open sets
are subsequently established
in \cite{CKS1, CKS2, CKS4, CKS5} for censored stable processes,
relativistic stable processes, mixed stable processes, and mixed Brownian motion and stable processes, respectively.
By integrating the two-sided heat kernel estimates in
Theorem \ref{t:main} with respect to $t$, we
obtain
the two-sided estimates on the Green function $ G_D(x, y):=\int_0^\infty
p_D(t, x, y)dt$ (see Theorem \ref{t:gfe} below), which extend
 \cite[Theorem 1.1]{KSV3}.

Condition {\bf (A)} is a very weak condition on the behavior of $\phi$
near infinity. Using the tables at the end of \cite{SSV},
 one can come up
plenty of explicit examples of  complete Bernstein functions satisfying
condition {\bf (A)}.
Here are a few of them.

(1) $\phi(\lambda)=\lambda^{\alpha/2}$, $\alpha\in (0, 2]$
(symmetric $\alpha$-stable process);

(2) $\phi(\lambda)=(\lambda+m^{2/\alpha})^{\alpha/2}-m$, $\alpha\in (0, 2)$ and $m>0$
(relativistic $\alpha$-stable process);

(3) $\phi(\lambda)=
\lambda^{\alpha/2}+ \lambda^{\beta/2} $, $0<\beta<\alpha <2$
(mixed symmetric $\alpha$- and $\beta$-stable processes);

(4) $\phi(\lambda)=\lambda^{\alpha/2}(\log(1+\lambda))^{p}$,
 $\alpha\in (0, 2)$, $p\in [-\alpha/2, (2-\alpha)/2]$.

\medskip

Now we give a way of constructing less explicit complete Bernstein functions
that have
very general
asymptotic behavior at infinity. Suppose that $\alpha\in (0, 2)$ and $\ell$ is a positive function
on $(0, \infty)$ which is slowly varying at infinity.
We further assume that $t \to t^{\alpha/2}\ell(t)$ is a right continuous increasing function with $\lim_{t \to 0} t^{\alpha/2}\ell(t)=0$ (so $\int^\infty_0(1+t)^{-2}t^{\alpha/2}\ell(t)dt < \infty$).
Then the function
$$
f(\lambda):=\int^\infty_0(\lambda+t)^{-2}t^{\alpha/2}\ell(t)dt
$$
is a Stieltjes function, and so the function
$$
\phi(\lambda):=\frac1{f(\lambda)}=\left(\int^\infty_0(\lambda+t)^{-2}t^{\alpha/2}\ell(t)dt\right)^{-1}
$$
is a complete Bernstein function (see \cite[Theorem 7.3]{SSV}). It follows from
\cite[Lemma 6.2]{WYY} that
$\phi (\lambda) \asymp \lambda^{\alpha /2} \ell (\lambda )$ when $\lambda\geq 2$.

The rest of the paper is organized as follows.
Section \ref{S:2} recalls and collects some preliminary
results that will be used in the sequel, including on-diagonal
heat kernel estimates and the boundary Harnack principle.
Section \ref{S:3} presents the interior lower bound heat kernel estimates,
including an explicit lower bound estimate for symmetric L\'evy
processes on $\R^d$. The proof of the short time factorization result
for $p_D(t, x, y)$ (that is, Theorem
\ref{thm:oppz}(i) and (ii))
is given in Section \ref{S:4}, while
the proof of
Theorem \ref{t:main}(i) and (ii) is given in Section \ref{S:5}.
The large time heat kernel estimates are proved in Section \ref{S:6}.
The Green function estimates for subordinate Brownian motions in bounded
$C^{1,1}$ open sets are derived in Section \ref{S:7} from the two-sided
Dirichlet heat kernel  estimates in Theorem \ref{t:main}.
The derivation, however, requires quite some effort.

\smallskip

Throughout this paper, $d\geq 1$
  and  the constants $R_0$, $R_1$, $R_2$, $\Lambda$, $\kappa$, $\delta_1$, $\delta_2$, $C_1$, $C_2$, $C_3$, $C_4$, $C_5$ and  $C_6$
  will be fixed.
 We use $c_1, c_2, \cdots$ to
denote generic constants, whose exact values are not important and
can  change from one appearance to another. The labeling of the
constants $c_1, c_2, \cdots$ starts anew in the statement of each
result. The dependence of the constant $c$ on the dimension $d$ will
not be mentioned explicitly. We will use ``$:=$" to denote a
definition, which is read as ``is defined to be". We will use
$\partial$ to denote a cemetery point and for every function $f$, we
extend its definition to $\partial$ by setting $f(\partial )=0$. We
will use $dx$ to denote the Lebesgue measure in $\bR^d$. For a Borel
set $A\subset \bR^d$, we also use $|A|$ to denote its Lebesgue
measure.

\section{Preliminary}\label{S:2}

In the first part of this section, we assume that $X$ is
a purely discontinuous rotationally symmetric L\'evy
process with L\'evy exponent $\Psi$ and L\'evy density
$J_X$ satisfying \eqref{e:psi2} and \eqref{e:psi1} respectively,
 where the
complete Bernstein function $\phi$ satisfies {\bf (A)}.

A function $u:\bR^d\mapsto [0, \infty)$ is said to be
harmonic in an open set $D\subset \bR^d$ with respect to $X$ if
for every open set $B$ whose closure is a compact subset of $D$,
\begin{equation}\label{e:har}
u(x)= \E_x \left[ u(X_{\tau_{B}})\right] \qquad
\hbox{for every } x\in B.
\end{equation}
A function $u:\bR^d\mapsto [0, \infty)$ is said to be regular
harmonic in an open set $D\subset \bR^d$ with respect to $X$ if
$$
u(x)= \E_x \left[ u(X_{\tau_{D}})\right]
\qquad \hbox{for every } x\in D.
$$
Clearly, a regular harmonic function in $D$ is harmonic in $D$.

Very recently the following form of the boundary Harnack principle is established in \cite{KSV5}.

\begin{thm}[{\cite[Theorem 1.1(i)]{KSV5}}]\label{ubhp}
There exists
 $c= c(\phi, \gamma )>0$ such that
for any $z_0 \in \R^d$, any open set $D\subset \R^d$, any $r\in (0,1)$
    and any nonnegative functions $u, v$ in $\R^d$ which are regular harmonic in $D\cap
    B(z_0, r)$ with respect to $X$ and vanish in $D^c \cap B(z_0, r)$, we have
$$
     \frac{u(x)}{v(x)}\,\le c\,\frac{u(y)}{v(y)}
     \qquad \mbox{ for all } x, y\in D\cap B(z_0, r/2).
$$
  \end{thm}

Recall also from \cite{KSV5} that $j$
enjoys the following properties:
for every $R>0$,
\begin{equation}\label{e:asmpbofjat0}
j(r)\asymp \frac{1}{r^d \, \Phi (r)}
\quad \hbox{for } r\in (0, R],
\end{equation}
and there exists $c>0$
such that
\begin{equation}\label{e:asmpbofjatinfty}
j(r)\le cj(r+1) \qquad \hbox{for } r\ge 1.
\end{equation}
It follows from \eqref{e:1.8}  and  \eqref{e:asmpbofjat0} that
both the function $\Phi$ defined by \eqref{e:1.9} and the function $j$
satisfy a doubling property; that is,
for every constant $R>0$, there is a constant
$c>1$ so that
\begin{equation}\label{e:j2}
\Phi (2r)\leq 4\, \Phi (r) \quad \hbox{and} \quad
j(r)\leq c j(2r) \quad \hbox{for every } r\in (0, R].
\end{equation}
Moreover, under condition {\bf (A)}, for any given $R>0$,
$j(r)$ satisfies all the conditions
in \cite{CK2} for $r\in (0, R)$. Thus
we have the following two-sided
estimates for $p(t, x, y)$ from
\cite{CKK2, CKK3}:

\begin{prop}\label{p:new1}
For any $T>0$, there exists $c_1=c_1(T, R, \gamma, \phi)>0$
such that
\begin{equation}\label{e:offdiag}
p(t, x, y) \leq
c_1 \, (\Phi^{-1}(t))^{-d}
\qquad  \hbox{ for } (t, x, y)\in (0, T]\times \R^d \times \R^d.
\end{equation}
For any $T, R>0$, there exists
$c_2=c_2(T, R, \gamma, \phi)>1$ such that
for all $(t, x, y)\in [0, T]\times \R^d\times \R^d$ with $|x-y| <R,$
\begin{align}\label{stssbound}
c^{-1}_2\left( (\Phi^{-1}(t))^{-d}\wedge tJ(x, y)\right)
\le p(t, x, y) \le
c_2 \left((\Phi^{-1}(t))^{-d}\wedge tJ(x, y)\right).
\end{align}
\end{prop}

\pf \eqref{e:offdiag} is given in the first display on page 1073 of \cite{CKK2}.
By \cite[Theorem 2.4]{CKK2} and \eqref{e:psi1}, there exist $T_*>0$ and $R_*$ such that
for all $(t, x, y)\in [0, T_*]\times \R^d\times \R^d$ with $|x-y| <R_*$,
\eqref{stssbound} holds.
Now we assume  $R>R_*$.
We construct $Z$ from $X$ by
removing jumps of size larger than $R$ via Meyer's construction (see \cite{Mey}).
Let $p_Z(t,x,y)$ be the transition density of $Z$.
By
\cite[Lemma 3.6]{BBCK} and \cite[Lemma 3.1(c)]{BGK} we have
for every $t>0$ and $x, y\in \R^d$,
\begin{equation}\label{e:offdiag1}
e^{-t\|\sJ_R\|_{\infty}}p_Z(t,x,y)\,\le\,  p(t,x,y)\, \le \, p_Z(t,x,y) +
 t \|J_R\|_{\infty},
\end{equation}
where
\begin{equation}\label{e:offdiag11}
J_R(x,y):=J_X(x,y){\bf 1}_{\{|x-y|> R\}} \quad \text{and}\quad
\sJ_{R}(x):=\int_{\R^d} J_R(x,y)dy.
\end{equation}
Applying \cite[Theorem 1.4]{CKK3} and its proof to \eqref{e:offdiag1}, we
 get
\begin{equation}\label{e:offdiag2}
c_1e^{-t\|\sJ_R\|_{\infty}}\left( (\Phi^{-1}(t))^{-d}\wedge tJ_X(x, y)\right)\,\le\,  p(t,x,y) \, \le \,c_2\left( (\Phi^{-1}(t))^{-d}\wedge tJ_X(x, y)\right) +
 t \|J_R\|_{\infty}.
\end{equation}
\eqref{stssbound} now follows from \eqref{e:psi1}. \qed

The function $J_X(x, y)$ gives rise to a L\'evy system for $X$, which
describes the jumps of the process $X$: for any non-negative
measurable function $f$ on  $\bR_+ \times \bR^d\times \bR^d$ with
$f(s, y, y)=0$ for all $y\in  \bR^d$
and stopping
time $T$ (with respect to the filtration of $X$),
\begin{equation}\label{e:levy}
\E_x \left[\sum_{s\le T} f(s,X_{s-}, X_s)
\right] = \E_x \left[ \int_0^T \left(
\int_{\bR^d} f(s,X_s, y) J_X(X_s,y) dy \right) ds \right].
\end{equation}
(See, for example, \cite[Proof of Lemma 4.7]{CK}
and \cite[Appendix A]{CK2}.)

When $\Psi(|\xi|)=\phi(|\xi|^2)$ and $\phi$ is a complete Bernstein function
satisfying condition {\bf (A)},
the following boundary Harnack principle on a $C^{1,1}$ open subset
with explicit decay rate is established in
\cite{KM} (see also \cite{KSV3}).

\begin{thm}[{\cite[Theorem 1.5]{KM}}]\label{t:bhp}
Suppose that $X$ is a rotationally symmetric L\'evy process with L\'evy exponent $\Psi(\xi)=
\phi(|\xi|^2)$ with $\phi$ being a complete Bernstein function satisfying
condition {\bf (A)}.
Assume that $D$ is a
(possibly unbounded) $C^{1, 1}$ open set in $\bR^d$ with
characteristics $(R_2, \Lambda)$. Then there exists
$c=c(R_2, \Lambda, \phi, d)>0$  such that for $r \in (0, (R_2 \wedge
1)/4]$, $Q\in \partial D$ and any nonnegative function $u$ in $\R^d$
that is harmonic in $D \cap B(Q, r)$ with respect to $X$ and
vanishes continuously on $ D^c \cap B(Q, r)$, we have
\begin{equation}\label{e:bhp_m}
\frac{u(x)}{\Phi(\delta_D(x))^{1/2}}\,\le c\,
\frac{u(y)}{\Phi(\delta_D(y))^{1/2}}
\qquad \hbox{for every } x, y\in  D \cap B(Q, r/2).
\end{equation}
\end{thm}

\section{Interior lower bound estimate}\label{S:3}

In this section, we assume that $X$ is
a purely discontinuous rotationally symmetric L\'evy
process with L\'evy exponent $\Psi$ and L\'evy density
$J_X$ satisfying \eqref{e:psi2} and \eqref{e:psi1} respectively,
 where the
complete Bernstein function $\phi$ satisfies {\bf (A)}.
We will give some preliminary lower bounds on $p_D(t, x, y)$
and $p(t, x, y)$.
We first recall the following from
 \cite{CKK2, CK2}.
Recall that $\sJ_R$ was defined in \eqref{e:offdiag11}.

\begin{lemma}\label{L:exit}
For any positive constant $M$,  there exist $\eps=\eps(M, \gamma, \phi)>0$ and $c=c(M, \gamma, \phi)\geq 1$
such that for all  $r \in (0, M]$ and $z\in \R^d$,
\begin{equation}\label{e:E1}
 \P_z \left( \tau_{B(z, r )} >
\eps \Phi (r) \right) \ge 2^{-1}
e^{-\|\sJ_1\|_{\infty}}
\end{equation}
and
\begin{equation}\label{e:E2}
 c^{-1}\, \Phi (r)\,\leq \,
 \E_z \left[ \tau_{B(z, r)} \right]\, \leq\, c \,\Phi (r).
\end{equation}
\end{lemma}

\pf  \eqref{e:E1} follows directly from
\cite[Lemma 2.5]{CKK2}.
We then have
$$
 \E_z \left[ \tau_{B(z, r)} \right]
\geq \eps \Phi (r) \P_z ( \tau_{B(z, r )} > \eps \, \Phi (r) ) >
2^{-1} e^{-\|\sJ_1\|_{\infty}}\eps \Phi (r) .
$$
On the other hand, by the L\'evy system for $X$ in \eqref{e:levy}, \eqref{e:psi1}, \eqref{e:asmpbofjat0}
and \eqref{e:j2},
we have
\begin{eqnarray*}
1&\geq & \P_z \left( X_{\tau_B(z, r)}\in B(z, 2r)^c \right)
=\E_z \left[ \int_0^{\tau_{B(z, r)}} \int_{B(z, 2r)^c} J_X(X_s, y) dy ds\right]\\
&\geq & c_1\E_z \left[ \int_0^{\tau_{B(z, r)}} \int_{B(z, 3r)\setminus B(z, 2r)} J(X_s, y) dy ds\right] \geq c_2 \frac{|B(z, 3r)\setminus B(z, 2r)|}{r^d \Phi (r)} \,
\E_z \left[ \tau_{B(z, r)} \right] ,
\end{eqnarray*}
which yields $\E_z \left[ \tau_{B(z, r)} \right]\leq c \Phi (r)$.
\qed

\begin{lemma}\label{l:cks2lemma3.1}
For any positive constants  $a$ and $R$, there exists
$c=c(a, R, \gamma, \phi)>0$ such that for all $z \in \bR^d$ and $r \in (0, R]$,
$$
\inf_{y\in B(z, r/2)}\P_y \left(\tau_{B(z,
r)}> a \Phi (r)  \right)\, \ge\, c.
$$
\end{lemma}

\pf  By \eqref{e:E1} and \eqref{e:j2}, there exists $\eps_1 =\eps_1 (R, \gamma, \phi)>0$ such that for all  $r\in (0, R]$,
$$
\inf_{z\in \bR^d } \P_z ( \tau_{B(z, r/2 )} >
\eps_1 \Phi (r) ) \ge
2^{-1}e^{-\|\sJ_1\|_{\infty}}.
$$
Thus it suffices to prove the lemma for $a>\eps_1$.  Applying the parabolic Harnack inequality \cite[Theorem 5.2]{CKK2}  at most
$2+ [ a/\eps_1 ]$ times, we conclude that there exists
$c_1=c_1(a, R, \gamma, \phi)>0$ such that
 for every  $w, y \in B(z,  r/2 )$,
$$
c_1\,p _{B(z, r )}(\eps_1 \Phi (r) ,z,w) \, \le
\, p _{B(z, r )}(a \Phi (r) ,y,w) .
$$
Thus
\begin{eqnarray*}
\P_y \left( \tau_{B(z, r)} > a \Phi (r) \right)
&=& \int_{B(z, r )}
p _{B(z, r)}(a \Phi (r) ,y,w) dw\\
&\ge& \int_{B(z, r/2 )} p _{B(z, r )}(a \Phi (r) ,y,w) dw \\
&\ge& c_1 \int_{B(z, r /2 )} p _{B(z, r /2 )}(\eps_1 \Phi (r) ,z,w) dw \\
& =& c_1 \P_z (\tau_{B(z, r /2)
}>\eps_1 \Phi (r))
\, \ge \, 2^{-1}e^{-\|\sJ_1\|_{\infty}} .
\end{eqnarray*}
This proves the lemma.
\qed

For the next four results, $D$ is an arbitrary nonempty open set and
we use the convention that $\delta_{D}(\cdot) \equiv \infty$ when $D
=\bR^d$.

\begin{prop}\label{p:cks2prop3.2}
Let $T>0$ and $a>0$   be constants.
  There exists
 $c=c(T,a,   \gamma, \phi)>0$ such that
 \bee\label{e:lb1}
p_D(t,x,y) \,\ge\,c\, (\Phi^{-1}(t))^{-d}
\eee
for every $(t, x, y)\in (0, T]\times D\times D$ with
$\delta_D (x) \wedge \delta_D (y)
 \ge a \Phi^{-1}(t) \geq 4 |x-y|$.
 \end{prop}

\pf
We fix
$(t, x, y)\in (0, T]\times D\times D$ satisfying
$\delta_D(x)\wedge \delta_D(y)\geq a \Phi^{-1} (t) \geq 4 |x-y|$.
Note that
$|x-y| \le a\Phi^{-1}(t)/4 \le a \Phi^{-1}(T)/4$
and that
$$
B(x, a\Phi^{-1}(t)/4) \subset B(y, a\Phi^{-1}(t)/2)\subset
  B(y, 2a\Phi^{-1}(t)/3)\subset D.
$$
So by the parabolic Harnack inequality \cite[Theorem 5.2]{CKK2}, there exists $c_1=c_1(T, \gamma, \phi)>0$ such that
$$
c_1 \, p_D(t/2, x, w) \, \le  \,  p_D(t,x,y) \quad \mbox{for every
} w \in B(x, a\Phi^{-1}(t)/4)  .
$$
This together with Lemma \ref{l:cks2lemma3.1} and \eqref{e:j2}
yields that \begin{eqnarray*}
p_D(t, x, y) &\geq & \frac{c_1}{ | B(x, a\Phi^{-1}(t)/4)|}
\int_{B(x, a\Phi^{-1}(t)/4)} p_D(t/2, x, w)dw\\
&\geq & c_2 (\Phi^{-1}(t))^{-d} \, \int_{B(x, a \Phi^{-1}(t)/4)}
p_{B(x, a\Phi^{-1}(t)/4)} (t/2, x, w)dw \\
&=& c_2 (\Phi^{-1}(t))^{-d} \, \P_x \left( \tau_{B(x, a\Phi^{-1}(t)/4)} >
t/2\right) \,\geq \, c_3 \, (\Phi^{-1}(t))^{-d},
\end{eqnarray*}
where $c_i=c_i(T, a, \gamma, \phi)>0$ for $i=2, 3$.
\qed

\begin{lemma}\label{l:cks2lemma3.3}
Let $T>0$ and $a>0$ be constants.
There exists
 $c=c ( a , T, \gamma, \phi )>0$ so that
$$
\P_x \left( X^{D}_t \in B \big( y, \,   a \Phi^{-1}(t)/2 \big)
\right) \,\ge \,c\,  (\Phi^{-1}(t))^d \, t J(x,y)
$$
for every $(t, x, y)\in (0, T]\times D\times D$ with $
\delta_D(x) \wedge \delta_D (y) \ge a \Phi^{-1}(t)$ and
$a \Phi^{-1}(t) \leq 4|x-y|$.

\end{lemma}

\pf
It follows from
Lemma \ref{l:cks2lemma3.1}
 that, starting at $z\in B(y, \,
a \Phi^{-1}(t)/4)$, with probability at least
$c_1=c_1
( a , T, \gamma, \phi )>0$
the process $X$ does not move more than $a \Phi^{-1}(t)/6$ by
time $t $. Thus, it suffices to show that there exists
$c_2=c_2(a, T, \gamma, \phi)>0$ such that for
$(t, x, y)\in (0, T]\times D\times D$ with
 $\delta_D(x) \wedge \delta_D (y) \ge a \Phi^{-1}(t)$ and $a \Phi^{-1}(t) \leq 4|x-y|$,
\begin{equation}\label{eq:molow}
\P_x \left(X^{D} \hbox{ hits the ball } B(y, \,
 a \Phi^{-1}(t)/4)\mbox{ by time } t \right)\, \ge\, c_2\,
(\Phi^{-1}(t))^d \, t   J(x,y) .
\end{equation}

Let $B^t_x:=B(x, \, a\Phi^{-1}(t)/9)$, $B^t_y:=B(y, \,  a \Phi^{-1}(t)/9)$
and $\tau^t_x:=\tau_{B^t_x}$. It follows from Lemma \ref{l:cks2lemma3.1}
 and \eqref{e:j2} that
there exists $c_3=c_3(a, T, \gamma, \phi )>0$ such that
\begin{equation}\label{eq:lowtau}
\E_x \left[ t \wedge
\tau^t_x \right] \,\ge\,t \, \P_x
\left(\tau^t_x \ge  t \right)
 \,\ge\, c_3\,t
\qquad \hbox{ for all } t\in (0, T].
\end{equation}
Since $B^t_x\cap B^t_y = \emptyset$, by the L\'evy system of $X$ and \eqref{e:psi1},
\begin{eqnarray}
&&\P_x \left(X^{D} \mbox{ hits the ball } B(y, \,
 a \Phi^{-1}(t)/4 )\mbox{ by time }   t \right) \nonumber\\
&\ge & \P_x\left(X_{t\wedge
\tau^t_x}\in B(y,\, a
\Phi^{-1}(t)/4) \right)
\ge  \E_x \left[\int_0^{t\wedge
\tau^t_x} \int_{B^t_y}
J_X (X_s, u) duds \right]\nn\\
&\ge & c_4\E_x \left[\int_0^{t\wedge \tau^t_x} \int_{B^t_y}
J (X_s, u) duds \right].
\label{e:dew}
\end{eqnarray}

We consider two cases separately.

(i) Suppose $|x-y| \le a \Phi^{-1}(T)$.
Since $|x-y| \ge   a \Phi^{-1}(t)/4$, we have for
 $s < \tau^t_x$ and $u \in B^t_y$,
$$
|X_s-u|\le |X_s-x|+|x-y|+|y-u| \le 2|x-y|.
$$
Thus by  \eqref{eq:lowtau} and \eqref{e:dew},
\begin{eqnarray*}
&&\P_x \left(X^{D} \mbox{ hits the ball } B(y, \,
  a \Phi^{-1}(t) )/4 \mbox{ by time }   t \right)\\
&\ge &c_4 \E_x \left[ t\wedge
\tau^t_x \right] \, |
B^t_y| \, j(2|x-y|)
\,\ge\, c_5\,  (\Phi^{-1}(t))^d \, t
j(2|x-y|)
\end{eqnarray*}
for some positive constant  $c_5=c_5(  a, T,  \gamma, \phi )>0$.
Therefore, in view of \eqref{e:j2},
the assertion of the lemma holds
when $|x-y| \le a \Phi^{-1}(T)$.

(ii) Suppose $|x-y| > a \Phi^{-1}(T)$.
In this case, for $s < \tau^t_x$ and $u \in B^t_y$,
\begin{eqnarray*}
&& |X_s-u|\, \le\,  |X_s-x|+|x-y|+|y-u|
\le |x-y| +
a \Phi^{-1}(t)/4 \le |x-y| +  a \Phi^{-1}(T)/4 .
\end{eqnarray*}
Thus from \eqref{e:dew} and then \eqref{eq:lowtau},
\begin{eqnarray*}
&&\P_x \left(X^{D} \mbox{ hits the ball } B(y, \,
a \Phi^{-1}(t)
/4)    \mbox{ by time }   t \right)\\
&\ge & c_6\E_x \left[ t\wedge
\tau^t_x \right] \int_{B^t_y} j\left(|x-y| +
a\Phi^{-1}(T)/4 \right)
\,  du\\
&\ge& c_7 \,  t\, | B^t_y| \,  j\left(|x-y| +
a\Phi^{-1}(T)/4 \right)
 \\
&\ge& c_8 \, t(\Phi^{-1}(t))^d j\left(|x-y| +
a\Phi^{-1}(T)/4 \right)
\end{eqnarray*}
for some constants
 $c_i=c_i ( a, T, \gamma, \phi )>0$, $i=6, 7, 8$.
Thus we conclude from
\eqref{e:asmpbofjatinfty} and  \eqref{e:j2} that the assertion of the lemma
holds for $|x-y| > a \Phi^{-1}(T)$ as well. \qed

\begin{prop}\label{p:cks2prop3.4}
Let $T$ and $a$ be positive constants. There exists $c=c(T, a, \gamma, \phi)>0$ such that
$$
p_D(t, x, y)\,\ge  \,c \,{ t}{J(x,y)}
$$
for every $(t, x, y)\in
(0, T]\times D\times D$ with $ \delta_D(x) \wedge \delta_D
(y) \ge a \Phi^{-1}(t)$ and
$a \Phi^{-1}(t) \leq 4|x-y|$.
\end{prop}

\pf
By the semigroup property,  Proposition \ref{p:cks2prop3.2},
 Lemma \ref{l:cks2lemma3.3} and \eqref{e:j2},
there exist positive constants
$c_1=c_1(T, a, \gamma, \phi)$ and $c_2=c_2(T, a, \gamma, \phi)$ such that
\begin{eqnarray*}
p_D(t, x, y) &=& \int_{D} p_{D}(t/2, x, z)
p_{D}(t/2, z, y)dz\\
&\ge& \int_{B(y, \, a \Phi^{-1}(t/2)/2)}
p_{D}(t/2, x, z) p_{D}(t/2, z, y) dz\\
&\ge& c_1 (\Phi^{-1}(t/2))^{-d} \,
\bP_x \left( X^{D}_{t/2} \in B(y, a
\Phi^{-1}(t/2)/2) \right)\\
&\ge & c_2\,{ t}{J(x,y)}.
\end{eqnarray*}
 \qed
\medskip

Combining Propositions \ref{p:cks2prop3.2} and \ref{p:cks2prop3.4}, we have the
following preliminary lower bound for $p_D(t,x,y)$.

\begin{prop}\label{step31}
Let $T$ and $a$ be
positive constants.
There exists $c=c(T, a, \gamma, \phi)>0$ such that
$$
p_D(t, x, y)\,\ge  \,c \, ((\Phi^{-1}(t))^{-d} \wedge  { t}{J(x,y)})
$$
for every $(t, x, y)\in
(0, T]\times D\times D$ with
$ \delta_D(x) \wedge \delta_D (y) \ge a \Phi^{-1}(t)$.
\end{prop}

In particular, Proposition \ref{step31} with $D=\R^d$ gives
\begin{thm}\label{t:globalhke}
For any constant $T>0$, there exists
$c=c(T, \gamma, \phi)>0$  such that
for all $(t, x, y)\in [0, T]\times \R^d\times \R^d$,
$$
p(t, x, y)\geq c\, \left((\Phi^{-1}(t))^{-d} \wedge tJ(x, y)\right).
$$
\end{thm}

\section{Factorization of Dirichlet heat kernel}
\label{S:4}

In this section,
we continue assuming  that $X$ is
a purely discontinuous rotationally symmetric L\'evy
process with L\'evy exponent $\Psi$ and L\'evy density
$J_X$ satisfying \eqref{e:psi2} and \eqref{e:psi1} respectively,
 where the
complete Bernstein function $\phi$ satisfies
the growth condition {\bf (A)}.

\begin{figure}[tb]
\begin{center}
\scalebox{0.60}{\includegraphics{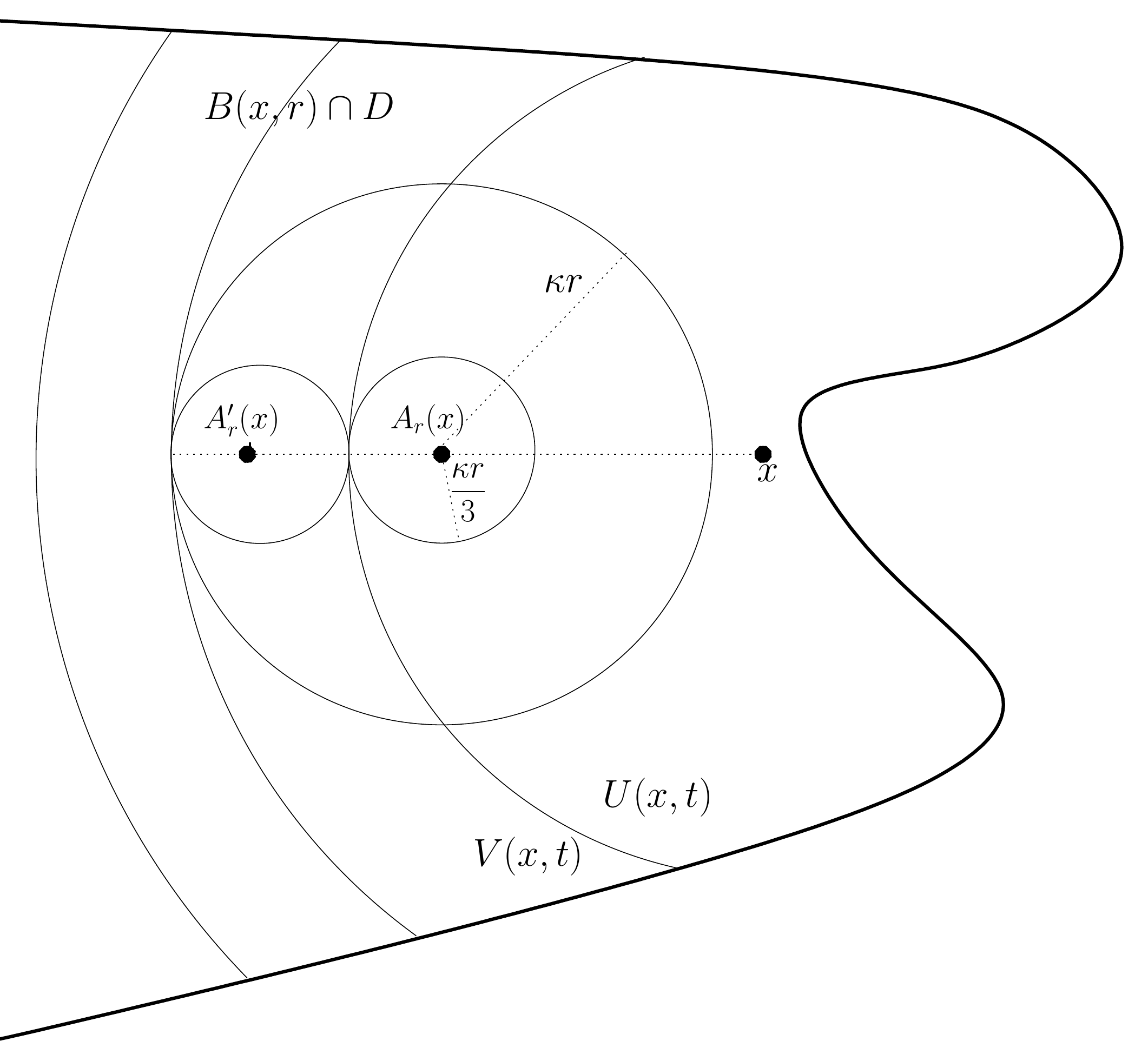}}
\caption{$U(x, t)$ and $V(x, t)$}
  \label{pic:ts}
  \end{center}
  \end{figure}

Throughout this section,
$T>0$ is a fixed constant and $D$ is a
fixed $\kappa$-fat open set with  characteristics $(R_1, \kappa)$. Recall
that $A_r(x)\in D$ is defined in Definition \ref{def:UB}.
For $(t, x) \in (0, T]\times \overline D$,
set $r=r(t)=\Phi^{-1}(t)R_1/\Phi^{-1}(T) \le R_1$
and define
\begin{equation}\label{e:UV}
 U(x, t):= D\cap B \big(x,|x-A_{r}(x)|+ \kappa r/3 \big), \quad
 V(x, t):= D\cap  B\big(x,|x-A_{r}(x)|+\kappa {r}\big).
\end{equation}
Let
   \begin{equation}\label{e:UVA}
A'_r(x)\in D \text{ be a point such that } B\big(A'_{r}(x), \kappa r/3 \big)
 \subset B \big( A_{r}(x),\kappa {r}\big)\setminus U(x, t).
 \end{equation}
Note that $ B(A_{r}(x), \kappa r/3) \subset U(x, t)$
and $ B(A'_{r}(x), \kappa \, r/ 6 )\subset B(A'_{r}(x), \kappa \, r/ 3 )
\subset V(x,t) \setminus U(x, t)$.
See Figure \ref{pic:ts}.

\begin{lemma}\label{lem:etd1r}
For every $T>0$ and $M \ge 1$, we have that, for $(t, x)\in (0, T]\times D$,
  \begin{eqnarray} \label{eq:cw}
  &&  \P_x(\tau_{D}>t/M) \asymp
        \P_x(\tau_{V(x, t)}>Mt) \asymp
     \P_x(\tau_{V(x, t)}>t/M) \asymp\P_x(\tau_{D}>Mt) \nonumber\\
  &  \asymp &
    \P_x(X_{\tau_{U(x, t)}}\in D) \asymp
   t^{-1} \E_x[\tau_{U(x, t)}],
  \end{eqnarray}
    where $U(x, t)$ and $V(x, t)$ are the sets defined in \eqref{e:UV} and
  the (implicit) comparison constants in  \eqref{eq:cw} depend
only on
$d, M, T, R_1, \kappa, \gamma$ and $\phi$.
\end{lemma}

\pf Without loss of generality we assume $R_1 \le 1$.
Fix
$(t, x)\in (0, T]\times D$, and set
$r=r(t)=\Phi^{-1}(t)R_1/\Phi^{-1}(T) \le R_1 \le 1$.
Recall that $U(x, t)$ and $V(x, t)$ are the sets defined in \eqref{e:UV}.
Observe that
\begin{equation}\label{e:v}
 \P_x (\tau_{V(x, t)}>Mt)\,\leq\, \P_x (\tau_{V(x, t)}>t/M) \wedge
\P_x (\tau_D>Mt) \,\leq\, \P_x (\tau_D>t/M).
\end{equation}
Note that by \eqref{e:psi1}, \eqref{e:asmpbofjat0} and (\ref{e:levy}), we have
 \begin{eqnarray}
&& \P_x\big(X({\tau_{U(x, t)}})\in B(A'_r(x), \kappa r/6)\big)
   \,=\,
   \E_x\int_0^{ \tau_{U(x, t)}} \int_{B(A'_r(x), \kappa r/6)}J_X(X_t,y)dtdy \nn\\
&&\asymp     r^d \, j(r)   \E_{ x}\big [\tau_{U(x, t)}\big ] \asymp     \, \Phi (r)^{-1} \, \E_{ x}\big [\tau_{U(x, t)}\big ]\,
 \asymp  \, t^{-1}\, \E_{ x}[\tau_{U(x, t)}] . \label{eq:cesp}
 \end{eqnarray}

 If $|x-A_{r}(x)|< \kappa r/2$,
 then $B(x, \kappa r/3)\subset U(x, t)\subset V(x, t)$ and so
\begin{eqnarray*}
1\,\geq \,
\P_x(\tau_D>t/M)\,\geq \,\P_x(\tau_D>Mt)
\,\geq\, \P_x (\tau_{V(x, t)}>Mt)
\,\geq\,
\P_x(\tau_{B(x, \kappa r/3 )}>Mt) ,
\end{eqnarray*}
which is greater than or equal to a positive constant depending only
on $\phi, T, R_1, \gamma,    d, \kappa, M$
by Lemma \ref{l:cks2lemma3.1}.
Thus we have, in view of \eqref{e:v},
$$
1 \asymp \P_x(\tau_{D}>t/M) \asymp
        \P_x(\tau_{V(x, t)}>Mt) \asymp
     \P_x(\tau_{V(x, t)}>t/M) \asymp
    \P_x(\tau_{D}>Mt).
$$
Moreover, when $|x-A_{r}(x)|< \kappa r/2$, we have
$B(x, \kappa r/3) \subset U(x, t)
\subset B(x, r)$ and so by
Lemma \ref{L:exit} and \eqref{e:j2},
$$c_1 t\leq
c_2 \Phi(\kappa r/3) \le \E_x[\tau_{B(x, \kappa r/3)}]
\leq \E_x[\tau_{U(x, t)}] \leq \E_x[\tau_{B(x,  r )}] \le
c_3\Phi(r)
\leq c_4 t.
$$
Combining the last two displays with \eqref{eq:cesp} and the fact that
$$
1\ge \P_x(X_{\tau_{U(x, t)}}\in D)\ge \P_x\big(
X_{\tau_{U(x, t)}}
\in B(A'_r(x), \kappa r/6)\big),
$$
we arrive at the assertion of the lemma when $|x-A_{r}(x)|<\kappa r /2 $.

Now we assume that $|x-A_{r}(x)| \ge \kappa r /2$.  We note that
\begin{equation}
  \label{eq:rsp}
  \P_x(\tau_{D}>t/M) \leq \P_x(\tau_{U(x, t)}>t/M)
  +\P_x(X_{\tau_{U(x, t)}}\in D)\,.
\end{equation}
For $r\in (0, R_1]$, by
Theorem \ref{ubhp},
we have
\begin{eqnarray*}
\P_x\big(X_{\tau_{U(x, t)}}\in D\big) &\leq& c_5
 \P_{A_r(x)}
\big(X_{\tau_{U(x, t)}}\in D\big)
\frac{\P_x\big(X_{\tau_{U(x, t)}}\in B(A'_r(x), \kappa r/6)\big)}
{\P_{A_r (x)}
\big(X_{\tau_{U(x, t)}}\in B(A'_r(x), \kappa r/6)\big)}\\
&\leq& c_5\frac{\P_x\big(X_{\tau_{U(x, t)}}\in B(A'_r(x), \kappa r/6)\big)}{
 \P_{A_r (x)}
\big(X_{\tau_{B(A_r(x), \kappa r/3)}}\in B(A'_r(x), \kappa r/6)\big)}.
\end{eqnarray*}
Note that when $(w, y) \in B(A_{r}(x),\kappa r/3)\times B(A'_r(x), \kappa r/6)$, $|w-y| \le 2 \kappa  r$
 and so
$j(|w-y|) \ge j(2 \kappa  r)$. Thus,
by \eqref{e:psi1}, \eqref{e:levy}
and Lemma \ref{L:exit},
\begin{eqnarray*}
&&\P_{A_{r}(x)} \big(X_{\tau_{B(A_r(x), \kappa r/3)}}\in B(A'_r(x), \kappa r/6)\big)
= \E_{A_{r}(x)}
\int_{B(A'_r(x), \kappa r/6)} \int_0^{\tau_{B(A_{r}(x),\kappa r/3)}} J_X(X_t,z)dtdz \\
&&\ge c_6  r^d \, j(r) \, \E_{A_{r}(x)}
\big[\tau_{B(A_{r}(x),\kappa r/3)}\big]
\geq c_7 \Phi (r)^{-1}\, \E_0\big[\tau_{B(0,\kappa r/3)}\big] \ge  c_8 \Phi (r)^{-1}\,\Phi(\kappa r/3) \ge c_9,
\end{eqnarray*}
where the last inequality is due to \eqref{e:j2}.
Therefore,
\begin{eqnarray}
\P_x\big(X_{\tau_{U(x, t)}}\in D\big)\,\leq\, c_{10}\, \P_x\big(X_{\tau_{U(x, t)}}\in B(A'_r(x), \kappa r/6)\big).\label{eq:cesp0}
\end{eqnarray}
Since
$\P_x\big (\tau_{U(x, t)}>t/M\big )\leq Mt^{-1} \E_x\big [\tau_{U(x, t)}\big ],$
we have by (\ref{eq:rsp}), \eqref{eq:cesp0}  and \eqref{eq:cesp},
\begin{eqnarray}
\P_x\big (\tau_{D}>t/M\big)\leq c_{11}  t^{-1}\E_x\big[\tau_{U(x, t)}\big] .\label{eq:cesp00} \end{eqnarray}  On the other hand, by the strong Markov property of $X$,
\begin{eqnarray*}
  &&\P_x\big ( X_{\tau_{U(x, t)}}\in B(A'_r(x), \kappa r/6)\big )
\\ &= & \E_x\left[ \P_{X_{\tau_{U(x, t)}}}\left(\tau_{B(X_{\tau_{U(x, t)}},
      \kappa r/6)}>Mt\right); X_{\tau_{U(x, t)}} \in B(A'_r(x), \kappa r/6)\right]\\
 &&+ \E_x\left[
    \P_{X_{\tau_{U(x, t)}}}\left(\tau_{B(X_{\tau_{U(x, t)}},\kappa r/6)}\le Mt \right): X_{\tau_{U(x, t)}} \in B(A'_r(x), \kappa r/6)\right]\\
    &\le & \E_x\left[
    \P_{X_{\tau_{U(x, t)}}}\left(\tau_{V(x, t)} >Mt\right); X_{\tau_{U(x, t)}} \in B(A'_r(x), \kappa r/6)\right]\\
    &&+\P_{0}\left(\tau_{B(0,\kappa r/6)}  \le Mt\right)
 \P_x\left( X_{\tau_{U(x, t)}} \in B(A'_r(x), \kappa r/6)\right)\\
  &= &
 \P_x\left(
    \tau_{V(x, t)} >Mt,  X_{\tau_{U(x, t)}} \in B(A'_r(x), \kappa r/6)\right)\\
    &&+\left(1-\P_{0}\left(\tau_{B(0,\kappa r/6)}>Mt\right)\right)
 \P_x\left( X_{\tau_{U(x, t)}} \in B(A'_r(x), \kappa r/6)\right)\\
   &\le  &
 \P_x\left(
    \tau_{V(x, t)} >Mt \right)+\left(1-\P_{0}\left(\tau_{B(0,\kappa r/6)}>Mt\right)\right)
 \P_x\left( X_{\tau_{U(x, t)}} \in B(A'_r(x), \kappa r/6)\right).
  \end{eqnarray*}
Thus
\begin{eqnarray}
\P_{0}\big(\tau_{B(0,\kappa r/6)}>Mt\big)
\P_x\left( X_{\tau_{U(x, t)}} \in B(A'_r(x), \kappa r/6)\right) \le \P_x\left(
    \tau_{V(x, t)} >Mt \right).\label{eq:cesp000}
\end{eqnarray}
Since $\P_{0}\big(\tau_{B(0,\kappa r/6)}>Mt\big)>c_{12}$ by
\eqref{e:j2} and Lemma \ref{l:cks2lemma3.1},
combining  \eqref{eq:cesp} with
\eqref{eq:cesp000}, we get
\begin{eqnarray*}
  t^{-1} \E_x[\tau_{U(x, t)}]
  \leq c_{13}\P_x(X_{\tau_{U(x, t)}}\in B(A'_r(x), \kappa r/6))
  \leq c_{14} \P_x(\tau_{V(x, t)}>Mt).
\end{eqnarray*}
This
together with \eqref{eq:cesp00} and \eqref{e:v} completes the proof
of this lemma.
\qed

\begin{lemma}\label{lemppu100}
 Suppose that $U_1,U_3, E$ are open subsets of $\bR^d$ with $U_1,
U_3\subset E$ and ${\rm dist}(U_1,U_3)>0$. Let $U_2 :=E\setminus
(U_1\cup U_3)$.  If $x\in U_1$ and $y \in U_3$, then for all $t >0$,
\begin{equation}\label{eq:ub}
p_{E}(t, x, y) \le \P_x\left(X_{\tau_{U_1}}\in U_2\right)
\left(\sup_{s<t,\, z\in U_2} p_E(s, z, y)\right)+ \left(t \wedge \E_x
\left[\tau_{U_1}\right] \right)\left(\sup_{u\in U_1,\, z\in U_3}
J_X(u,z)\right)
\end{equation}
  and
  \begin{equation}\label{eq:lb}
p_{E}(t, x, y)\,\ge\, t\,  \P_x(\tau_{U_1}>t)
\,\P_y(\tau_{U_3}>t)\inf_{u\in U_1,\, z\in U_3}
J_X(u,z).
\end{equation}
\end{lemma}
\pf
For \eqref{eq:ub}, see the proof of \cite[Lemma 3.4]{CKS5}.
For \eqref{eq:lb}, see the proof of \cite[Lemma 3.3]{CKS6}.
\qed

\noindent{\bf Proof of Theorem \ref{thm:oppz}(ii).}
It follows from Lemma \ref{lem:etd1r} and the semigroup property that
it suffices to establish the assertion for $T\le 1$. So in the
remainder of this proof, we assume that $T\le 1$.
Fix $t\in (0, T]$ and set $r=\Phi^{-1}(t) R_1/\Phi^{-1}(T)$.
By \eqref{e:asmpbofjat0}, \eqref{e:j2}, \eqref{e:offdiag} and Theorem \ref{t:globalhke},
   \begin{eqnarray}
   p(t/2,x,y)\asymp (\Phi^{-1}(t))^{-d} \quad \text{if }|x-y|\leq 8  {r}.
    \label{eq:ppun22} \end{eqnarray}
By the semigroup property,  \eqref{e:j2}, \eqref{e:offdiag}
and Lemma~\ref{lem:etd1r},  when $|x-y| \le 8  {r}$,
  \begin{eqnarray}
    p_{D}(t/2,x,y)&=&\int_{D}
 p_D(t/4,x,z)p_D(t/4,z,y)dz \nn\\
& \leq&
        \sup_{z \in \R^d}
     p(t/4,z,y) \P_x(\tau_{D}>t/4)
    \nonumber\\
   &\leq&
   c_1 (\Phi^{-1}(t))^{-d}\P_x(\tau_{D}>t).\label{eq:ppun222}
  \end{eqnarray}
Thus by \eqref{eq:ppun22} and \eqref{eq:ppun222},
\begin{eqnarray}
    p_{D}(t/2,x,y)
    \le
    c_2\,\P_x(\tau_{D}>t)p(t/2,x,y)\quad \text{if }|x-y|\leq 8  {r}\,.\label{eq:ppun2}
  \end{eqnarray}

  Now we assume $|x-y|>8  {r}$.
  Let
  $D_1:=U(x, t)$ be the set defined
  in \eqref{e:UV},
  $D_3:= \{z\in D: |z-x|>|x-y|/2\}$ and
  $$D_2 :=D\setminus
(D_1\cup D_3)=\{z\in D \setminus U(x, t): |z-x| \le |x-y|/2\}.$$
 Then by  condition {\bf (B)} we have
 $$
 \sup_{s<t/2,\, z\in D_2} p(s, z, y)
  \le   \sup_{s<t/2,\,   |z-y| \geq |x-y|/2} p(s, |z-y|)\\
  \le     C_1 \sup_{s<t/2} p(s, C_2 |x-y|/2).
$$
 Extend the definition of $p(t, r)$ by setting $p(t, r)=0$ for $t<0$
 and $r\geq 0$. For each fixed $x, y\in \R^d$ and $t>0$ with $|x-y|> 8r$,
 one can easily check, in view of \eqref{e:1.7}, that $(s, w)\mapsto p(s, C_2 |w-y|/2)$
is a parabolic function in $(-\infty, T]\times B(x, 2r)$.
So by the parabolic Harnack inequality from \cite[Theorem 5.2]{CKK2},
 there is a constant
$c_3=c_3(d,  \gamma,  \phi)\geq 1$
 so that for every $t\in (0, T]$,
$$
\sup_{s<t/2} p(s, C_2 |x-y|/2)\leq  c_3 p(t/2, C_2 |x-y|/2).
$$
Hence we have
\begin{equation}\label{eq:ppun31}
 \sup_{s<t/2,\, z\in D_2} p(s, z, y) \leq
 c_3  C_1 p(t/2, C_2|x-y|/2)
 \leq c_3 C_1^2 p(t/2, C_2^2 |x-y|/4),
 \end{equation}
where in the last inequality, we applied condition {\bf (B)} again.
If $u\in D_1$, then $|x-u| \le  |x-A_r(x)|+\kappa r/3$. Thus
     if $u\in D_1,\, z\in
    D_3$, then
  \begin{eqnarray*}
  |u-z|  \ge |z-x|-|x-u| \ge |z-x|-|x-A_r(x)|-\kappa r/3  \ge
  |z-x|-{r} \ge \frac{1}{2}|z-x| \ge \frac{1}{4}|x-y|.
   \end{eqnarray*}
Thus  by \eqref{e:asmpbofjat0}, \eqref{e:psi1}, \eqref{e:j2} and Theorem \ref{t:globalhke},
\begin{eqnarray*}
   t\sup_{u\in D_1,\, z\in D_3}J_X(u,z)
& \le&  \gamma  t \sup_{ |u-z| \ge |x-y|/4}j(|u-z|)
 \leq   \gamma t  j(|x-y|/4)  \leq  \gamma t  j(C_2^2 |x-y|/4)\\
&\leq& c_4 p(t/2, C_2^2 |x-y|/4).
\end{eqnarray*}
 Let $C_5:=C_2^2/4$.
 It follows then from \eqref{eq:ppun31}, Lemmas \ref{lem:etd1r} and \ref{lemppu100}
 that for  $|x-y|>8  {r}$,
  \begin{eqnarray}
    p_{D}(t/2, x, y)&\le&
    c_5p(t/2,C_5 x, C_5 y)\left(\P_x(X_{\tau_{U(x, t)}}\in D)
    + t^{-1}  \E_x[\tau_{U(x, t)}]\right)\nn\\
    &\leq& c_6\P_x(\tau_{D}>t)\,p(t/2, C_5x, C_5y) . \label{eq:ppun3}
  \end{eqnarray}

   Hence by condition {\bf (B)}, \eqref{eq:ppun2}, \eqref{eq:ppun3}, symmetry, the semigroup
  property and Lemma~\ref{lem:etd1r}, we conclude that for all $(t, x, y)\in (0, T]\times D\times D$,
  \begin{eqnarray*}
    p_D(t,x,y)&=&\int_D p_D(t/2,x,z)p_D(t/2,z,y)dz\\
    &\leq&
    c_6^2\P_x(\tau_{D}>t)\P_y(\tau_{D}>t)\int_{\R^d} p(t/2,C_5x,C_5z)p(t/2,C_5z,C_5y)dz\\
     &=&   c_{7}\P_x(\tau_{D}>t)\P_y(\tau_{D}>t)\int_{\R^d} p(t/2,C_5x,z)p(t/2,z,C_5y)dz\\
    &= &c_{7}\P_x(\tau_{D}>t)\P_y(\tau_{D}>t)p(t,C_5x,C_5y)\,.
  \end{eqnarray*}
\qed

\noindent{\bf Proof of Theorem \ref{thm:oppz}(i).}
Fix $(t, x, y)\in (0, T]\times D \times D$
and set $r=\Phi^{-1}(t)R_1/\Phi^{-1}(T)$.
Let $U(x, t)$ be the set defined
in \eqref{e:UV},  $A_r(x)\in D$ and $A'_r(x)\in D$ be the points
defined in Definition \ref{def:UB} and \eqref{e:UVA} respectively.
 Note that
 \begin{eqnarray}
  p_D(t, x, y)&=& \int_{D\times D} p_D(t/3, x, u) p_D(t/3, u, v) p_D (t/3, u, y) du dv \nn \\
  &\geq &  \inf_{(u,v) \in B(A'_r(x), \kappa r/6) \times
B(A'_r(y), \kappa r/6)} p_D(t/3, u, v) \int_{B(A'_r(x), \kappa r/6)}   p_{D}(t/3,x,u)du \, \nn\\
&& \cdot  \int_{B(A'_r(y), \kappa r/6)}p_D(t/3,v,y)dv. \label{e:pD}
\end{eqnarray}
For $(u,v) \in B(A'_r(x), \kappa r/6) \times B(A'_r(y), \kappa r/6)$,
\begin{eqnarray}
&& \inf_{(u,v) \in B(A'_r(x), \kappa r/6) \times
B(A'_r(y), \kappa r/6)} p_D(t/3, u, v) \nn \\
&\geq & c_1 \inf_{(u,v) \in B(A'_r(x), \kappa r/6) \times
B(A'_r(y), \kappa r/6)} \left( (\Phi^{-1}(t/3))^{-d} \wedge (t/3) J(u, v) \right) \nn \\
&\ge&  c_2 (tJ(x,y)\wedge (\Phi^{-1}(t))^{-d}), \label{e:pD2}
\end{eqnarray}
where in the first inequality we used Proposition \ref{step31}, in the second inequality we
used \eqref{e:asmpbofjat0} and \eqref{e:asmpbofjatinfty} respectively in the cases
$|x-y| \ge \kappa r$ and $|x-y| < \kappa r$.

On the other hand, for $u\in B(A'_r(x), \kappa r/6)$, by Lemma~\ref{lemppu100} with
  $U_1=U(x, t)$ and
  $U_3=B(A'_{r}(x),\kappa r/4 )$,
  we  have
 \begin{eqnarray*}
p_{D}(t/3,x,u)
&\geq& t\P_x(\tau_{U(x, t)}>t/3) \P_u (\tau_{U_3}>t/3) \inf_{w\in
  U(x, t),\,z\in U_3}
  J_X(w,z) \\
&\geq & \gamma^{-1} t\P_x(\tau_{U(x, t)}>t/3) \P_u (\tau_{B(u, \kappa r /16)}>t/3)
\inf_{w\in U(x, t),\,z\in U_3}j(|w-z|) .
\end{eqnarray*}
Since
$$  j(|w-z|) \ge j(r)=j(\Phi^{-1}(t)R_1/\Phi^{-1}(T))
\quad \hbox{for } (w,z)\in U(x, t) \times U_3,
$$
we have by Lemma \ref{l:cks2lemma3.1}, \eqref{e:asmpbofjat0}
and \eqref{e:j2},
 \begin{eqnarray*}
p_{D}(t/3,x,u)
  \geq c_3 t\, \P_x(\tau_{U(x, t)}>t/3)  \frac1{  \Phi^{-1}(t)^d \, t  }= c_3  \Phi^{-1}(t)^{-d}  \, \P_x(\tau_{U(x, t)}>t/3)  .
  \end{eqnarray*}
It follows from condition {\bf (A)} that there exists $c_4>1$ such that
$\Phi^{-1}(bt) \le c_4 b^{1/(2\delta_2)}\Phi^{-1}(t)$ for every $t \le T$ and $b \in (0, 1]$.
 Thus with $a:=(\kappa/(3c_4))^{2\delta_2}$ we have that
 $3^{-1} \kappa \Phi^{-1}(t) \ge \Phi^{-1}(a t)$ for every $t \le T$,
and hence
$V(x,at) \subset U(x,t)$.
Thus Lemma \ref{lem:etd1r} implies that
$$
\P_x(\tau_{U(x, t)}>t/3) \ge  \P_x(\tau_{V(x, at)}>t/3)  \ge c_{5} \P_x(\tau_{D}>t ).
$$
Consequently, by \eqref{e:j2},
$$
\int_{B(A'_r(x), \kappa r/6)}   p_{D}(t/3,x,u)du
\geq \frac{c_{6}}{\Phi^{-1}(t)^d}  \P_x(\tau_{D}>t ) | B(A'_r(x), \kappa r/6)|
\geq c_{7} \P_x(\tau_D>t ).
$$
Similarly, using the symmetry we also have
 $\int_{B(A'_r(y), \kappa r/6)}p_D(t/3,v,y)dv\geq c_{8} \P_y(\tau_D>t )$.
Hence we have by \eqref{e:pD} and \eqref{e:pD2},
\begin{eqnarray*}
    p_{D}(t,x,y)\geq c_{9}\, \P_x(\tau_{D}>t) \P_y(\tau_{D}>t)
  \left((\Phi^{-1}(t))^{-d}\wedge tJ(x,y)\right) .
    \end{eqnarray*}
\qed

\section{Small time Dirichlet heat kernel estimates for subordinate Brownian motions in $C^{1,1}$ open set}\label{S:5}

In this section we assume that $D$ is a $C^{1,1}$ open set in $\bR^d$ with
characteristics $(R_2, \Lambda)$
and  that $X$ is
a subordinate Brownian motion
 with L\'evy exponent $\Psi(\xi)=
\phi(|\xi|^2)$, where $\phi$ is a complete Bernstein function satisfying
condition {\bf (A)}.

\bigskip

\noindent{\bf Proof of Theorem \ref{t:main}(ii).}
Without loss of generality we assume $R_2 \le 1$.
In view of Theorem \ref{thm:oppz}(ii), it suffices to show that
\begin{equation}\label{e:5.1}
\P_{x}\left(\tau_{D}>t\right) \le  c \left( 1\wedge \frac{\Phi(\delta_D(x))}{t} \right)^{1/2}
\quad \hbox{for } (t, x) \in (0, T]\times D.
\end{equation}
Fix $(t, x) \in (0, T]\times D$ and set
$r=r(t)= \Phi^{-1}(t)R_2/\Phi^{-1}(T) \le R_2 \le 1$.
By Theorem \ref{t:globalhke}, we only need to show the theorem for
 $\delta_{D}(x)<r/16$. Take $x_0\in \partial D$
such that $\delta_{D}(x)=|x-x_0|$.

Let
$U_1:=B( x_0, r/8 ) \cap D$ and
${\bf n}(x_0)$ be the unit inward normal of $\partial{D}$ at
the point $x_0$. Put
$$x_1=x_0+\frac{r}{16 }{\bf n}(x_0).$$
 Note that
$\delta_{D} (x_1)=r/16$.
Applying the boundary Harnack principle (Theorem \ref{t:bhp}) and \eqref{e:j2} we get
\begin{align}
 \P_x(X_{\tau_{U_1}}\in D \setminus U_1) &\le c_1
  \P_{x_1}(X_{\tau_{U_1}}\in D \setminus U_1)  \sqrt{\frac{\Phi(\delta_D(x))} { \Phi(\delta_D(x_1))}}\nn\\
 & \le  c_2
  \P_{x_1}(X_{\tau_{U_1}}\in D \setminus U_1)
  \sqrt{\frac{\Phi(\delta_D(x))}{ t}}
  \le  c_2
  \sqrt{\frac{ \Phi(\delta_D(x))} {t}}. \label{e:nnw11}
\end{align}
Take $x_2\in \R^d$ so that
$B(x_2, r) \subset B(x_0 , 4r) \setminus B(x_0 ,r).$
Then, by \eqref{e:levy},
\eqref{e:asmpbofjat0} and \eqref{e:j2}, we have
 \begin{eqnarray*}
\P_x ( X_{\tau_{U_1}} \in B(x_2, r)  )
& =&   \E_x \left[
\int_0^{\tau_{U_1} } \int_{B(x_2, r) } J (|X_s-y|) dy ds\right]\\
&\ge& c_3
|B(0, \Phi^{-1}(t))| j(5r)
  \E_x [\tau_{U_1}]\\
  &\ge& c_4  \Phi(\Phi^{-1}(t))^{-1}
  \E_x [\tau_{U_1}] = c_4  t^{-1}
  \E_x [\tau_{U_1}].
 \end{eqnarray*}
Now by the same argument as that in \eqref{e:nnw11}, we get,
 \begin{eqnarray}
  \E_x [\tau_{U_1}] \le   c_4^{-1}   t \P_x ( X_{\tau_{U_1}} \in B(x_2, r)  ) \le
c_5t  \P_{x_1}(X_{\tau_{U_1}} \in B(x_2, r)) \sqrt{\frac{ \Phi(\delta_D(x))}
 {t}}
\le c_5\sqrt{t\, \Phi(\delta_D(x))} . \label{e:nnw2}
 \end{eqnarray}
Thus, by \eqref{e:nnw11} and \eqref{e:nnw2}, we have
\begin{eqnarray}
\P_{x}\left(\tau_{D}>t\right) &\le& \P_x\left(\tau_{U_1} >t\right) +
\P_x \left( X_{\tau_{U_1}} \in D\setminus U_1\right) \nn\\
&\le& \frac1{t} \E_x\left[\tau_{U_1}\right] +\P_x \left( X_{\tau_{U_1}} \in D\setminus U_1\right)
\le c_6 \sqrt{\frac{ \Phi(\delta_D(x))}{t}}.
\end{eqnarray}
Since $\P_{x}\left(\tau_{D}>t\right)\le 1$, \eqref{e:5.1} follows immediately.
\qed

Let $\delta_{\partial D}(x)$ be  the Euclidean
distance between $x$ and $\partial D$.
 It is well-known that any $C^{1, 1}$ open set $D$
 with $C^{1,1}$-characteristics $(R_2, \Lambda)$
satisfies both the {\it uniform interior ball condition} and the
{\it uniform exterior ball condition}: there exists
$r_0=r_0(R_2, \Lambda)\leq R_2$
such that for every $x\in D$ with $\delta_{\partial D}(x)< r_0$ and $y\in
\bR^d \setminus \overline D$ with $\delta_{\partial D}(y)<r_0$,
there are $z_x, z_y\in \partial D$ so that $|x-z_x|=\delta_{\partial
D}(x)$, $|y-z_y|=\delta_{\partial D}(y)$ and that $B(x_0,
r_0)\subset D$ and $B(y_0, r_0)\subset \bR^d \setminus \overline D$
for $x_0=z_x+r_0(x-z_x)/|x-z_x|$ and $y_0=z_y+r_0(y-z_y)/|y-z_y|$.

In the remainder of this section, we fix such an $r_0$ and set
$T_0:=\Phi( r_0/16)$.
For any $x \in D$ with $\delta_D(x) < r_0$,
let $z_x$ be a point on $\partial D$  such that $|z_x-x|=\delta_D (x)$
and ${\bf n}(z_x):=(x-z_x)/|z_x-x|$.

\begin{lemma}\label{l:cks2lemma4.6}
Let $ \kappa_0 \in (0, 1)$ and $a>0$.
There exists a constant $c=c(
\kappa_0, R_2, r_0, a,  \phi)>0$
such that for every  $(t, x)\in (0, T_0]\times D$
with $\delta_D(x) \leq 3 \Phi^{-1}(t) < r_0/4$ and $\kappa_0 \in (0, 1)$,
\begin{equation}\label{e:4.5}
\P_x \left( X^{D}_{at} \in B(x_0, \kappa_0 \Phi^{-1}(t)) \right)
\,\ge \, c\, \sqrt{\frac{ \Phi(\delta_D(x))}{t}},
\end{equation}
where  $x_0:=z_x+ \tfrac{9}{2}\Phi^{-1}(t){\bf n}(z_x)$.
\end{lemma}

\pf Let $0< \kappa_1 \le  \kappa_0$ and assume first that
$2^{-4}\kappa_1 \Phi^{-1}(t) < \delta_D(x)\leq 3\Phi^{-1}(t)$. As in
the proof of Lemma \ref{l:cks2lemma3.3}, we get that, in this case, using
the fact that $|x-x_0|\in [
 \tfrac{3}{2} \kappa_0 \Phi^{-1}(t), 6\Phi^{-1}(t)]$,
there exist constants $c_i=c_i(\kappa_1, r_0, a)>0$,
$i=1,2,$ such that for all  $t\le T_0$, we have
\begin{equation}\label{e:case1}
\P_x \left( X^{D}_{at} \in B(x_0, \kappa_1 \Phi^{-1}(t)) \right)
\,\ge \,c_1 t(\Phi^{-1}(t))^dJ(x, x_0) \,\ge \,c_2>0 .
\end{equation}
By taking
$\kappa_1= \kappa_0$, this shows that \eqref{e:4.5} holds for all $a>0$
in the case when
$2^{-4} \kappa_0 \Phi^{-1}(t) < \delta_D(x)\leq 3\Phi^{-1}(t)$.

So it suffices to consider the case that $\delta_D(x) \leq
2^{-4} \kappa_0 \Phi^{-1}(t)$. We now show that there is some $a_0>1$
so that \eqref{e:4.5} holds for every $a\geq a_0$ and $\delta_D(x)
\leq 2^{-4}
\kappa_0 \Phi^{-1}(t)$. For simplicity, we assume without
loss of generality that $x_0=0$ and let $\wh B:=B(0, \kappa_0
\Phi^{-1}(t))$.  Let $U:=D\cap B(z_x,  \kappa_0 \Phi^{-1}(t))$.
By the strong Markov property of $X^{D}$ at the first exit time
$\tau_{U}$ from $U$ and Lemma \ref{l:cks2lemma3.1}, there exists
$c_3=c_3(a)>0$ such that
\begin{eqnarray}
&& \P_{x} \left(X^{D}_{at} \in {\wh B} \right) \nonumber \\
&\geq & \P_{x}\left( \tau_{U}<at, \ X_{\tau_{U}} \in B(0,
2^{-1} \kappa_0 \Phi^{-1}(t)) \hbox{ and }
|X^D_s-X_{\tau_{U}}|<2^{-1} \kappa_0 \Phi^{-1}(t) \hbox{ for }
s\in [\tau_U, \tau_U+at]\right) \nonumber \\
&\geq& c_3 \, \P_{x}\left( \tau_U<at \hbox{ and }
X_{\tau_U}\in B(0, 2^{-1} \kappa_0 \Phi^{-1}(t)) \right)
\label{e:4.8}.
\end{eqnarray}
Let  $x_1
=z_x+4^{-1} \kappa_0 {\bf n}(z_x)\Phi^{-1}(t)$
and  $B_1:= B(x_1, 4^{-1} \kappa_0 \Phi^{-1}(t))$.
It follows from  the boundary Harnack principle (Theorem \ref{t:bhp})
and \eqref{e:j2}
that there exist
$c_k= c_k( R_2, \Lambda, \gamma, \phi)>0$, $k=4, 5$,
such that for all $t\in (0, T_0]$,
\begin{eqnarray*}
\P_{x} \left(X_{\tau_U}\in B(0, 2^{-1} \kappa_0
\Phi^{-1}(t))\right)
&\ge& c_4
\P_{x_1} \left(X_{\tau_U}\in B(0,
2^{-1} \kappa_0 \Phi^{-1}(t))\right)
\sqrt{\frac{\Phi(\delta_D(x))}{\Phi(\delta_D(x_1))}}\\
&\ge& c_5\P_{x_1} \left(X_{\tau_{B_1}}\in B(0,
2^{-1} \kappa_0 \Phi^{-1}(t))\right)
\sqrt{\frac{\Phi(\delta_D(x))}{t}}.
\end{eqnarray*}
By \eqref{e:levy}, Lemma \ref{L:exit} and
\eqref{e:asmpbofjat0}-\eqref{e:j2}, we have
\begin{eqnarray*}
&& \P_{x_1} \left(X_{\tau_{B_1}}\in B(0, \kappa_0 \Phi^{-1}(t)/2) \right)
= \E_x \left[ \int_0^{\tau_{B_1}} \int_{B(0, \kappa_0 \Phi^{-1}(t)/2)}
J(X_s, y) dy \right] \\
&\geq & c_6 j(\kappa_0 \Phi^{-1}(t)) \, |B(0, \kappa_0 \Phi^{-1}(t) /2)| \,
\E_x \left[ \tau_{B_1} \right] \\
& \geq&
\frac{c_7}{ ( \kappa_0 \Phi^{-1}(t))^d \Phi ( \kappa_0 \Phi^{-1}(t))}
 \, ( \kappa_0 \Phi^{-1}(t))^d \, \Phi (\kappa_0 \Phi^{-1}(t)) =c_7.
\end{eqnarray*}
Thus
\begin{equation}\label{e:case4}
\P_{x} \left(X_{\tau_U}\in B(0, 2^{-1} \kappa_0
\Phi^{-1}(t))\right)
\ge c_8 \sqrt{\frac{\Phi(\delta_D(x))}{t}}.
\end{equation}
It follows from \eqref{e:nnw2} that there exists $
c_9>0$ such that
\begin{eqnarray*}
\P_{x}(  \tau_U\geq a t ) \,\leq \, (at)^{-1} \, \E_{x} [
\tau_U ] \le \,
a^{-1} c_9\, \sqrt{\frac{\Phi(\delta_D(x))}{t}}.
\end{eqnarray*}
Define $a_0= 2c_9/(c_8) $. We  have by
 \eqref{e:4.8}--\eqref{e:case4}
and the display above that for $a\geq a_0$ ,
\begin{eqnarray}
\P_x ( X^D_{at}\in \wh B ) &\geq& c_3 \, \left( \P_{x}  (
X_{\tau_U}\in B(0, 2^{-1}\kappa_0 \Phi^{-1}(t) ) ) -  \P_{x}
\left(  \tau_U\geq at  \right) \right)  \nonumber\\
&\geq&  c_3 \, (c_9/2) \,
\sqrt{\frac{\Phi(\delta_D(x))}{t}}.\label{e:case2}
\end{eqnarray}
\eqref{e:case1} and \eqref{e:case2} show that \eqref{e:4.5} holds
for every $a\geq a_0$ and for every $x\in D$ with $\delta_D(x) \leq
3 \Phi^{-1}(t)$.

Now we deal with the case $0<a<a_0$ and $\delta_D(x) \leq
2^{-4} \kappa_0 \Phi^{-1}(t)$. If $\delta_D(x) \leq 3 \Phi^{-1}(at/a_0)$,
we have from \eqref{e:4.5} for the case of
$a=a_0$ that there exist
$c_{10}=c_{10}( \kappa_0, R_2, \Lambda, a, \phi)>0$ and
$c_{11}=c_{11}( \kappa_0, R_2, \Lambda, a, \phi)>0$ such that
\begin{eqnarray*}
\P_x \left( X^D_{at}\in B(x_0,  \, \kappa_0 \Phi^{-1}(t))
\right) &\ge & \P_x \left( X^D_{a_0 (at/a_0)}\in B(x_0,   \,
\kappa_0  \Phi^{-1}(at/a_0)) \right) \\
&\geq&  c_{10} \, \sqrt{\frac{\Phi(\delta_D(x))}{at/a_0}}= c_{11} \,
\sqrt{\frac{\Phi(\delta_D(x))}{t}}.
\end{eqnarray*}
If $3 \Phi^{-1}(at/a_0) < \delta_D(x) \leq 2^{-4} \kappa_0
\Phi^{-1}(t)$ (in this case $1>\kappa_0 > 3 \cdot 2^4
\Phi^{-1}(a/a_0)$), we get \eqref{e:4.5} from \eqref{e:case1} by
taking $\kappa_1=\Phi^{-1}(a/a_0)$. The proof of the lemma is now
complete. \qed

\medskip

\noindent {\bf Proof of  Theorem \ref{t:main}(i)}.
By Theorem \ref{thm:oppz}(i), it suffices to show that
\begin{equation}\label{e:5.10}
\P_{x}\left(\tau_{D}>t\right) \geq  c \left( 1\wedge \frac{\Phi(\delta_D(x))}{t} \right)^{1/2} \quad \hbox{for }
(t, x) \in (0, T]\times D.
\end{equation}
Assume  $(t, x) \in (0, T]\times D$.
Since $D$ satisfies the uniform interior ball condition
with radius $r_0$ and $0< {(T_0/T)t} \le T_0$,
we can choose  a point $\xi^t_x$ as follows:
if $\delta_D(x) \le 3 \Phi^{-1}(
{(T_0/T)t})$,
 let $\xi^t_x=z_x+ (9/2)\Phi^{-1}(
{(T_0/T)t}){\bf
n}(z_x)$ so that
$$B(\xi^t_x, (3/2)  \Phi^{-1}(
{(T_0/T)t}))\subset B(
z_x+ 3\Phi^{-1}(
{(T_0/T)t}){\bf
n}(z_x), 3\Phi^{-1}(
{(T_0/T)t})) \setminus \{x\}$$ and $\delta_D(z) \ge 3\Phi^{-1}(
{(T_0/T)t})$
for every $z \in B(\xi^t_x, (3/2)  \Phi^{-1}(
{(T_0/T)t}))$).
If $\delta_D(x) > 3 \Phi^{-1}(
{(T_0/T)t})$, choose $\xi^t_x \in B(x, \delta_D(x))$
so that $|x-\xi^t_x|=(3/2)  \Phi^{-1}(
{(T_0/T)t})$. Note that in this case,
$$B(\xi^t_x, (3/2) \Phi^{-1}(
{(T_0/T)t}))\subset B(x, \delta_D(x)) \setminus \{x\}$$ and $\delta_D(z) \ge \Phi^{-1}( {(T_0/T)t})$
for every $z \in B(\xi^t_x,  2^{-1} \Phi^{-1}({(T_0/T)t}))$.
We also define $\xi^t_y$ the same way.

If $\delta_D(x) \le 3 \Phi^{-1}( {(T_0/T)t})$, by
Lemma \ref{lem:etd1r} (with $M=T/T_0$ when $T \ge T_0$) and Lemma \ref{l:cks2lemma4.6} (with $a= 1, \kappa=2^{-1}$),
\begin{align*}
\P_x\left( \tau_D >
{t}\right) \ge
c_1\P_x\left( \tau_D > {(T_0/T)t}\right) \ge c_1
\P_x\left(X^D_{ {(T_0/T)t}} \in B(\xi^t_x, 2^{-1}
\Phi^{-1}( {(T_0/T)t}) )\right)  \ge c_2 \sqrt{\frac{\Phi(\delta_D(x))}{t}}.
\end{align*}
If  $\delta_D(x) >  3 \Phi^{-1}({(T_0/T)t})$,
by Lemma \ref{lem:etd1r}, Proposition \ref{step31} and
\eqref{e:j2},
\begin{eqnarray}
&&  \P_x\left( \tau_D >
{t}\right) \ge c_1\P_x\left(X^D_{
{(T_0/T)t}} \in B(\xi^t_x, 2^{-1}
\Phi^{-1}({(T_0/T)t}) )\right)  \nonumber \\ &=& c_1\int_{B(\xi^t_x, 2^{-1}\Phi^{-1}({(T_0/T)t}) )} p_{{D}}({(T_0/T)t},x,u)du
\,\ge\, c_3  \,\ge\, c_4 \left(1\wedge \frac{\Phi(\delta_D(x))}{t} \right)^{1/2}.    \nn
\end{eqnarray}
\qed

\section{Large time heat kernel estimates}\label{S:6}

In this section, we
first give
 the proofs of Theorems
\ref{thm:oppz}(iii) and \ref{t:main}(iii).

\bigskip

\noindent {\bf Proof of Theorem \ref{thm:oppz}(iii).} \
Since $D$ is bounded, in view of \eqref{e:offdiag},
the transition semigroup $\{P^D_t, t>0\}$
of $X^D$ consists of Hilbert-Schmidt operators, and hence compact operators,
in $L^2(D; dx)$.
So $P^D_t$ has discrete spectrum
$\{e^{-\lambda_k t}; k\geq 1\}$, arranged in decreasing order and repeated
according to their multiplicity.
Let $\{\phi_k, k\geq 1\}$ be the corresponding eigenfunctions with
unit $L^2$-norm ($\|\phi_1\|_{L^2(D)}=1$)
which forms an orthonormal basis for $L^2(D; dx)$.

Clearly, for every $k \ge 1$
\begin{equation}\label{ne:4.2}
\int_D \P_x(\tau_D > 1)
\phi_k (x) dx   \le |D|^{1/2} \|\phi_k\|_{L^2(D)} =|D|^{1/2}.
\end{equation}
By using the eigenfunction expansion of $p_D$ we get
\begin{equation}\label{ne:4.3}
\int_{D\times D} \P_x(\tau_D > 1) p_D(t, x, y)  \P_y(\tau_D > 1) \, dx dy =
\sum_{k=1}^\infty e^{-t \lambda_k} \left( \int_D \P_x(\tau_D > 1)
\phi_k (x) dx\right)^2.
\end{equation}
 Noting that
$\lambda_k$ is increasing and $
\| f\|_{L^2( D)}^2
= \sum_{k=1}^\infty (\int_D f(z)\phi_k(z)dz)^2$,
 we have
for all $t>0$,
\begin{align}
\int_{D\times D} \P_x(\tau_D > 1) p_D(t, x, y)  \P_y(\tau_D > 1) \, dx dy
&\leq
e^{-t \lambda_1} \, \int_D \P_x(\tau_D > 1)^2 dx \nn\\
&\leq
 e^{-t \lambda_1} \, |D|.\label{ne:4.4}
\end{align}

On the other hand,
by Theorem \ref{thm:oppz}(ii), Remark \ref{R:1.1}(iii)
and \eqref{ne:4.2} we have that
there is a constant $c_1>0 $ so that for every  $x\in D$,
\begin{align}
&\phi_1 (x)= e^{\lambda_1}
\int_D p_D (1, x, y)
\phi_1 (y) dy\nn\\
&\leq c_1 \P_x(\tau_D > 1)\int_D \P_y(\tau_D > 1)
\phi_1 (y) dy \leq   c_1 |D|^{1/2}\, \P_x(\tau_D > 1) .
\end{align}
It now follows from \eqref{ne:4.3}  that for every that for every
 $t>0$,
\begin{eqnarray}\label{ne:4.55}
&&  \int_{D\times D} \P_x(\tau_D > 1) p_D(t, x, y)  \P_y(\tau_D > 1) \, dx dy\geq
e^{-t \lambda_1} \, \left(\int_D  \P_x(\tau_D > 1) \phi_1(x) dx
\right)^2 \nonumber  \\
&&\geq  e^{-t \lambda_1} \, \left(\int_D  c_1^{-1}    |D|^{-1/2}
\phi_1(x)^2 dx\right)^2 = c_1^{-2}    |D|^{-1} \, e^{-t \lambda_1} .
\end{eqnarray}

For $t\geq 3$ and $x, y\in D$, we have  that
\begin{equation}\label{ne:4.6}
p_D (t, x, y)=\int_{D\times D}
p_D (1, x, z) p_D(t-2, z, w)
 p_D(1, w, y) dz dw .
\end{equation}
By Theorem \ref{thm:oppz}(ii), Remark \ref{R:1.1}(iii),
\eqref{e:offdiag}
and \eqref{ne:4.4} we have that there are constants
$c_i>0$, $i=2, 3$, so that for every
$t\geq 3$ and $x, y\in D$,
\begin{align}
p_D(t, x, y) &\leq c_2\P_x(\tau_D > 1)\P_y(\tau_D > 1)  \int_{D\times D}
 \P_z(\tau_D > 1) p_D (t-2, z, w) \P_w(\tau_D > 1) dz dw \nn\\
 &\leq   c_3 \, \P_x(\tau_D > 1)\P_y(\tau_D > 1) e^{-t\lambda_1}. \label{ne:4.7}
\end{align}
By  Theorem \ref{thm:oppz}(i),
Theorem \ref{t:globalhke},
the boundedness of $D$ and \eqref{ne:4.55}  we have
that there are constants
$c_i>0$, $i=4, 5$,
so that for every  $t\geq 3$ and $x, y\in D$,
\begin{align*}
p_D(t, x, y) &\geq  c_4\, \P_x(\tau_D > 1)\P_y(\tau_D > 1)
 \int_{D\times D} \P_z(\tau_D > 1) p_D (t-2, z, w) \P_w(\tau_D > 1) dz dw\\
& \geq c_5 \, \P_x(\tau_D > 1)\P_y(\tau_D > 1) e^{-t\lambda_1}.
\end{align*}
This combined with \eqref{ne:4.7} establishes Theorem \ref{thm:oppz}(iii). \qed

\noindent {\bf Proof of Theorem \ref{t:main}(iii).} \
By Theorem \ref{t:main}(i),
it suffices to prove the theorem  for $T\geq
3$.
By \eqref{e:5.1}, \eqref{e:5.10} and the boundedness of $D$, $\P_x(\tau_D > 1) \asymp  \Phi(\delta_D(x))^{1/2}$. This and Theorem \ref{thm:oppz}(iii) imply Theorem \ref{t:main}(iii).
\qed

\section{Green function estimates}\label{S:7}

In this section,
we use Theorem \ref{t:main} to get sharp two-sided estimates on the
 Green functions of subordinate Brownian motions in bounded $C^{1, 1}$ open sets.
We first establish the following two lemmas.
\begin{lemma}\label{l:phi_a}
For every $r\in (0, 1]$ and every open subset $U$ of $\R^d$,
\begin{align}
\frac12\left(1\wedge \frac{r^2\Phi (\delta_U(x))^{1/2}\Phi(\delta_U(y))^{1/2}}{\Phi(|x-y|)}\right)
&\leq
\left(1\wedge \frac{r \Phi (\delta_U(x))^{1/2}}{\Phi(|x-y|)^{1/2}}\right)
\left(1\wedge \frac{r \Phi (\delta_U(y))^{1/2}}{\Phi(|x-y|)^{1/2}}\right)\nonumber\\
&\leq    1\wedge \frac{r^2\Phi (\delta_U(x))^{1/2}\Phi(\delta_U(y))^{1/2}}{\Phi(|x-y|)} .\label{e:12}
\end{align}
\end{lemma}

\pf  The second inequality holds trivially.
Without loss of generality, we assume $\delta_U(x)\leq \delta_U(y)$.
If both $\frac{r \Phi (\delta_U(x))^{1/2}}{\Phi(|x-y|)^{1/2}}$
and $\frac{r \Phi (\delta_U(y))^{1/2}}{\Phi(|x-y|)^{1/2}}$ are less than
1 or if both are
larger than one,
$$
\left(1\wedge \frac{r \Phi (\delta_U(x))^{1/2}}{\Phi(|x-y|)^{1/2}}\right)
\left(1\wedge \frac{r \Phi (\delta_U(y))^{1/2}}{\Phi(|x-y|)^{1/2}}\right)
=  1\wedge \frac{r^2\Phi (\delta_U(x))^{1/2}\Phi(\delta_U(y))^{1/2}}{\Phi(|x-y|)} .
$$
So we only need to consider the case when $\frac{r \Phi
(\delta_U(x))^{1/2}}{\Phi(|x-y|)^{1/2}}\leq 1< \frac{r \Phi
(\delta_U(y))^{1/2}}{\Phi(|x-y|)^{1/2}}$.
 Note that
$ \Phi (\delta_U (y)) \leq \Phi (\delta_U(x)+|x-y|)$.
If $\delta_U(x)\geq |x-y|$, then by
\eqref{e:1.8},
$\Phi (\delta_U(y))\leq \Phi (2 \delta_U (x)) \leq 4  \Phi (\delta_U(x)) $
and so
$$  1\wedge \frac{r^2\Phi (\delta_U(x))^{1/2}\Phi(\delta_U(y))^{1/2}}{\Phi(|x-y|)}
\leq 1\wedge
 \frac{ 2 r \Phi
(\delta_U(x)) }{\Phi(|x-y|) }
 \leq 2 \left( 1\wedge \frac{r \Phi (\delta_U(x))^{1/2}}{\Phi(|x-y|)^{1/2}}\right).
 $$
When $\delta_U(x)<|x-y|$,  by \eqref{e:1.8} again, $\Phi (\delta_U(y))
\leq \Phi (2 |x-y|) \leq 4  \Phi(|x-y|) $ and so
$$
1\wedge \frac{r^2\Phi (\delta_U(x))^{1/2}\Phi(\delta_U(y))^{1/2}}{\Phi(|x-y|)}
\leq 1\wedge \frac{2 r^2\Phi (\delta_U(x))^{1/2}  }{\Phi(|x-y|)^{1/2}}
\leq 2 \left( 1\wedge \frac{r \Phi (\delta_U(x))^{1/2}}{\Phi(|x-y|)^{1/2}}\right),
$$
where the assumption $r \le 1$ is used in the last inequality.
This establishes the first inequality of \eqref{e:12}. \qed

By condition {\bf (A)}, we have that for every $T>0$,
 there exist $C_T>1$ such that
\begin{equation} \label{e:77}
C_T^{-1}\left(\frac{r}{R} \right)^{1/(2 \delta_1)} \le \frac{ \Phi^{-1} ( r)}{\Phi^{-1} ( R)} \le C_T\left(\frac{r}{R} \right)^{1/(2 \delta_2)}
 \quad \hbox{for }
0<r \le R\le T.
\end{equation}
Moreover, for every $M>0$, we have
\bee\label{e:6.20}
 r\Phi '(r) \asymp \Phi (r) \quad \hbox{for } r\in (0, M]
 \eee
(see the paragraph after \cite[Lemma 1.3]{KM}).

\begin{lemma}\label{l:new1dim}
Suppose $T>0$ and set
\bee\label{e:7.4}
h_T(a,r)=a+ \Phi(r) \int_{\Phi(r)/T}^1 \left(1 \wedge \frac{ua}{\Phi(r)}\right) \frac{1}{u^2\Phi^{-1}(u^{-1} \Phi(r))}\, du+\frac{\Phi(r)}{r}\left(1 \wedge \frac{a}{\Phi(r)}\right).
\eee
Then
$$ h_T(a, r) \asymp \frac{a}{r} \wedge \left( \frac{a}{\Phi^{-1}(a)}
+ \left( \int_{r}^{\Phi^{-1}(a)} \frac{\Phi(s)}{s^2}ds \right)^+ \, \right)
$$
for  $0<r\leq \Phi^{-1}(T/2)$  and $0<a \le  (2^{-1} \wedge (2C_T)^{-2\delta_2}) T$,
 where $C_T$ is the constant in \eqref{e:77} and $x^+:=x\vee 0$.
\end{lemma}

\pf
 For $(a, r)$ with
 $0<a <  \Phi(r)\le   T/2 $,
\begin{align*}
h_T (a,r) \asymp
a+ a\int_{\Phi(r)/T}^1  \frac{du}{u\Phi^{-1}(u^{-1} \Phi(r))}
  + \frac{a}{r}
  =a+ \frac{a}{r}\int_{\Phi(r)/T}^1  \frac{\Phi^{-1}(\Phi(r))}{\Phi^{-1}(u^{-1} \Phi(r))}u^{-1}du + \frac{a}{r}.
\end{align*}
By \eqref{e:77}, since $\Phi(r) \le T/2$, we have
\begin{align*}
0<c_2= c_1^{-1}\int^{1}_{1/2}
 u^{\frac{1}{2\delta_1}-1} du  \le \int_{\Phi(r)/T}^1  \frac{\Phi^{-1}(\Phi(r))}{\Phi^{-1}(u^{-1} \Phi(r))}u^{-1}du
 \le  c_1\int^{1}_{0}
 u^{\frac{1}{2\delta_2}-1} du =c_3<\infty.
\end{align*}
Thus,
for  $0<a <  \Phi(r)\le   T/2 $, we have
\bee \label{e:7.5}
c_2\frac{a}{r} \le h_T(a,r)
\le c_3 \left( a+  \frac{a}{r}\right) \le  c_4 \frac{a}{r}.
\eee

On the other hand,
for $(a, r)$ with
$\Phi(r)\leq a \le (2^{-1} \wedge (2C_T)^{-2\delta_2}) T$,
using the change of variable
$u=\Phi (r)/\Phi (s)$ and then applying integration by parts
for the first integral below,  we have
 \begin{eqnarray}
 h_T(a,r) & \asymp &
a+ \Phi(r) \int_{ \Phi(r)/a}^1 \frac{du}{u^2\Phi^{-1}(u^{-1} \Phi(r))}+a\int_{\Phi(r)/T}^{\Phi(r)/a} \frac{du}{u\Phi^{-1}(u^{-1} \Phi(r))}
+\frac{\Phi (r)}{r}
\nn\\
&=& a+ \int_r^{\Phi^{-1}(a)} \frac{\Phi'(s)}{s} ds+ a\int_{\Phi^{-1}(a)}^{\Phi^{-1}(T)}
\frac{\Phi '(s)}{s\Phi (s)} ds + \frac{\Phi (r)}{r} \nn \\
&=& a+ \left(\frac{a}{\Phi^{-1}(a)}-\frac{\Phi(r)}{r}\right)
+ \int_{r}^{\Phi^{-1}(a)} \frac{\Phi(s)}{s^2}ds +
  a \int_{\Phi^{-1}(a)}^{\Phi^{-1}(T)}
\frac{\Phi '(s)}{s\Phi (s)} ds + \frac{\Phi (r)}{r} \nn \\
&=&  a+ \frac{a}{\Phi^{-1}(a)} + \int_{r}^{\Phi^{-1}(a)} \frac{\Phi(s)}{s^2}ds
 +  a \int_{\Phi^{-1}(a)}^{\Phi^{-1}(T)}
\frac{\Phi '(s)}{s\Phi (s)} ds  .
 \label{e:gt1}
\end{eqnarray}
Since
$a\le \left(2^{-1}\wedge (2C_T)^{-2\delta_2}\right)T$,
by \eqref{e:77} and the fact that $\Phi^{-1}$ is increasing,
\begin{align}
 \frac{1}{\Phi^{-1}(a)}-\frac{1}{\Phi^{-1}(T)}\asymp \frac{1}{\Phi^{-1}(a)}\ge c_4\label{e:gt2}
\end{align}
for some $c_4>0$.
Using  \eqref{e:6.20} and  \eqref{e:gt2} in the second integral
   in \eqref{e:gt1}, we get that for $(a, r)$ with
$\Phi(r)\leq a \le (2^{-1} \wedge (2C_T)^{-2\delta_2}) T$,
\begin{eqnarray}
  h_T (a,r)&\asymp&  a+  \frac{a}{\Phi^{-1}(a)}
+ \int_{r}^{\Phi^{-1}(a)} \frac{\Phi(s)}{s^2}ds+
a\int_{\Phi^{-1}(a)}^{\Phi^{-1}(T)}
\frac{1}{s^2} ds  \nn \\
&=& a+  \frac{a}{\Phi^{-1}(a)}  +\int_{r}^{\Phi^{-1}(a)} \frac{\Phi(s)}{s^2}ds+a
 \left( \frac{1}{\Phi^{-1}(a)}-\frac{1}{\Phi^{-1}(T)} \right) \nn \\
& \asymp&  \frac{ a}{\Phi^{-1}(a)}
+\int_{r}^{\Phi^{-1}(a)} \frac{\Phi(s)}{s^2}ds.  \label{e:7.8}
\end{eqnarray}

Since $\Phi (s)$ is an increasing function,
when $0<\Phi (r) \leq a$, we have
$$\frac{ a}{\Phi^{-1}(a)}
+\int_{r}^{\Phi^{-1}(a)} \frac{\Phi(s)}{s^2}ds
\leq \frac{ a}{\Phi^{-1}(a)}
+ a \int_{r}^{\Phi^{-1}(a)} \frac{1}{s^2}ds = \frac{a}{r},
$$
while when $\Phi (r)\geq a>0$,
$$\frac{ a}{\Phi^{-1}(a)}
+ \left( \int_{r}^{\Phi^{-1}(a)} \frac{\Phi(s)}{s^2}ds \right)^+
= \frac{ a}{\Phi^{-1}(a)} \geq   \frac{a}{r} .
$$
This combined with \eqref{e:7.5} and \eqref{e:7.8} establishes the lemma.
\qed

Recall that the Green function $G_D(x, y)$ of $X$ in $D$ is defined as $G_D(x, y)=\int^\infty_0p_D(t, x, y)dt$.

\begin{thm}\label{t:gfe}
Suppose that $X$ is a subordinate Brownian motion with L\'evy exponent $\Psi(\xi)=
\phi(|\xi|^2)$ with $\phi$ being a complete Bernstein function satisfying
condition {\bf (A)}.
Let $D$ be a bounded $C^{1,1}$ open subset of $\bR^d$
with characteristics $(R_2, \Lambda)$
and $a(x,y)=\Phi (\delta_D(x))^{1/2}
  \Phi (\delta_D(y))^{1/2}$, $x,y \in D$.

\begin{description}
\item{\rm (i)}
There exists
$c_1 >0$ depending only on $\text{diam}(D),   R_2, \Lambda, d$ and  $\phi$
such that for all $d\ge 1$ and $(x, y)\in
D\times D$,
$$
G_D(x, y)\ge c_1 \frac{\Phi(|x-y|)} {|x-y|^{d}}
 \left(1\wedge \frac{  \Phi (\delta_D(x))}{ \Phi(|x-y|)}\right)^{1/2}
 \left(1\wedge \frac{  \Phi (\delta_D(y))}{ \Phi(|x-y|)}\right)^{1/2} .
$$

\item{\rm (ii)}
There exists
$c_2 >0$ depending only on $\text{diam}(D), R_2, \Lambda,  d$ and  $\phi$
such that for
all $d\ge 1$ and $(x, y)\in
D\times D$,
$$
G_D(x, y)\le c_2 \frac{ a(x,y)}{ |x-y|^d} .
$$
\item{\rm (iii)}
Let $d=1$. Then for
 $(x, y)\in D\times D$,
$$
G_D(x, y)\asymp
  \frac{ a(x,y)}{ |x-y|} \wedge \left( \frac{a(x,y)}
  {\Phi^{-1}(a(x,y) )}+ \left( \int_{|x-y|}^{\Phi^{-1}(a(x,y))} \frac{\Phi(s)}{s^2}ds \right)^+ \, \right).
$$

\item{\rm (iv)}
Let $d\geq 2$. Then for $(x, y)\in D\times D$,
\begin{eqnarray*}
G_D(x, y)&\asymp& \frac{\Phi(|x-y|)} {|x-y|^{d}}
 \left(1\wedge \frac{  \Phi (\delta_D(x))}{ \Phi(|x-y|)}\right)^{1/2}
 \left(1\wedge \frac{  \Phi (\delta_D(y))}{ \Phi(|x-y|)}\right)^{1/2}
\\&\asymp& \frac{\Phi(|x-y|)} {|x-y|^{d}} \left(1\wedge \frac{a(x,y)}{ \Phi(|x-y|)} \right).
\end{eqnarray*}
\end{description}
\end{thm}

\pf
    Put $T=
        (2 \vee (2C_T)^{2\delta_2}) \Phi({\rm diam} (D))$,
where $C_T$ is the constant in \eqref{e:77}.
It follows from Theorem \ref{t:main}(iii) that
\begin{equation}\label{e:newgfe}
\int^\infty_Tp_D(t, x, y)dt \asymp a(x,y).
\end{equation}

Using the boundedness of $D$, Remark \ref{R:1.1}(iii) and \eqref{e:asmpbofjat0},
the results of Theorem \ref{t:main}(i)--(ii)
can be rewritten as follows: there exists $c_1>0$ such that for $(t, x, y)\in (0, T]\times D\times D$,
\begin{eqnarray}
&&c_1^{-1}\left(1\wedge \frac{\Phi (\delta_D(x))}t \right)^{1/2} \left(1\wedge \frac{\Phi (\delta_D(y))}t \right)^{1/2}
\left( (\Phi^{-1}(t))^{-d}\wedge \frac{t}{|x-y|^d\Phi(|x-y|)}\right)\nn\\
&\leq & p_D(t, x, y) \label{e:nlo} \\
&\leq &c_1\left(1\wedge \frac{\Phi (\delta_D(x))}t \right)^{1/2} \left(1\wedge \frac{\Phi (\delta_D(y))}t \right)^{1/2}
\left( (\Phi^{-1}(t))^{-d}\wedge \frac{t}{|x-y|^d\Phi(|x-y|)}\right)\label{e:nup} .
\end{eqnarray}

By the change of variable
$u= \frac{\Phi(|x-y|)}{t}$ and the fact that $t\to \Phi^{-1}(t)$ is increasing, we have
\begin{eqnarray}
&& \int_0^{T} \left(1\wedge \frac{\Phi (\delta_D(x))}t \right)^{1/2} \left(1\wedge \frac{\Phi (\delta_D(y))}t \right)^{1/2}
\left( (\Phi^{-1}(t))^{-d}\wedge \frac{t}{|x-y|^d\Phi(|x-y|)}\right)  dt
\nonumber\\
&=&\frac{\Phi(|x-y|)} {|x-y|^{d}}  \left(\int_{\Phi(|x-y|)/T}^1+\int_1^\infty\right)
 u^{-2}\left(   \left( \frac{ \Phi^{-1} (  ut)}{   \Phi^{-1} (  t)}  \right)^d
\wedge u^{-1} \right)  \left(1\wedge \frac{ {\sqrt u} \Phi (\delta_D(x))^{1/2} }{ \Phi(|x-y|)^{1/2} }\right) \nonumber\\
&&\qquad \qquad \qquad \qquad \qquad \qquad \qquad \times\left(1\wedge \frac{
{\sqrt u} \Phi (\delta_D(y))^{1/2}} {\Phi(|x-y|)^{1/2} }\right) du
\nonumber\\
&\asymp&
 \frac{\Phi(|x-y|)} {|x-y|^{d}}   \int_{\Phi(|x-y|)/T}^1  u^{-2}\left(  \frac{ |x-y|}{   \Phi^{-1} (  u^{-1} \Phi(|x-y|))} \right)^d
  \left(1\wedge \frac{  u a(x,y) }{ \Phi(|x-y|) }\right)  du
 \nonumber\\
&& + \frac{\Phi(|x-y|)} {|x-y|^{d}} \int_1^\infty
 u^{-3}  \left(1\wedge \frac{ {\sqrt u} \Phi (\delta_D(x))^{1/2} }{ \Phi(|x-y|)^{1/2} }\right) \left(1\wedge \frac{
{\sqrt u} \Phi (\delta_D(y))^{1/2}} {\Phi(|x-y|)^{1/2} }\right) du
\nonumber\\
&=:&I+II,        \label{e:2}
\end{eqnarray}
where in the fourth line of the display above,
we used Lemma \ref{l:phi_a}.

\medskip

(i) The estimate on $II$ is easy.
\begin{align}
&\frac1{2 }  \frac{\Phi(|x-y|)} {|x-y|^{d}}  \left(1\wedge \frac{ \Phi (\delta_D(x))^{1/2}}{ \Phi(|x-y|)^{1/2} }\right) \left(1\wedge \frac{ \Phi (\delta_D(y))^{1/2}}{ \Phi(|x-y|)^{1/2} }\right)\nonumber\\
&= \frac{\Phi(|x-y|)} {|x-y|^{d}}  \int_{1}^\infty
  u^{-3} \, \left(1\wedge
\frac{   \Phi (\delta_D(x))^{1/2}}{ \Phi(|x-y|)^{1/2} }\right)
\left(1\wedge \frac{  \Phi (\delta_D(y))^{1/2}}{
\Phi(|x-y|)^{1/2} }\right) du \nonumber\\
 \le \, & II  \,\le \,
 \frac{\Phi(|x-y|)} {|x-y|^{d}}  \int_{1}^\infty
  u^{-2} \, \left(u^{-1/2} \wedge
\frac{   \Phi (\delta_D(x))^{1/2}}{ \Phi(|x-y|)^{1/2} }\right)
\left(u^{-1/2} \wedge \frac{  \Phi (\delta_D(y))^{1/2}}{
\Phi(|x-y|)^{1/2} }\right) du \nonumber\\
&\leq
 \frac{\Phi(|x-y|)} {|x-y|^{d}}
\int_{1}^\infty  u^{-2} \,
  \left(1 \wedge
\frac{   \Phi (\delta_D(x))^{1/2}}{ \Phi(|x-y|)^{1/2} }\right) \left(1
\wedge \frac{  \Phi (\delta_D(y))^{1/2}}{
\Phi(|x-y|)^{1/2} }\right) du \nonumber\\
& = \frac{\Phi(|x-y|)} {|x-y|^{d}}    \left(1\wedge \frac{
\Phi (\delta_D(x))^{1/2}}{ \Phi(|x-y|)^{1/2} }\right) \left(1\wedge \frac{
\Phi (\delta_D(y))^{1/2}}{ \Phi(|x-y|)^{1/2} }\right)  . \label{e:4}
\end{align}
Now part (i) of the theorem follows from  \eqref{e:nlo} and the lower bound of $II$ in \eqref{e:4}.
\medskip

(ii) We let
\bee \label{e:6.14}
u_0:= \frac{ a(x,y)}{ \Phi(|x-y|)}.
\eee
Clearly $1/u_0 \geq \Phi (|x-y|)/\Phi ({\rm diam}(D))\ge 2\Phi (|x-y|)/T$.
By \eqref{e:77},
\begin{eqnarray}
 I &\le & \frac{ a(x,y)}{|x-y|^d}\int^{1}_{\Phi(|x-y|)/T}
 \frac{ |x-y|^d }{   \Phi^{-1} (  u^{-1} \Phi(|x-y|))^d}\, u^{-1} du \nn\\
 &= & \frac{ a(x,y)}{|x-y|^d}\int^{1}_{\Phi(|x-y|)/T}
 \left(\frac{ \Phi^{-1} (  \Phi(|x-y|))}{   \Phi^{-1} (  u^{-1} \Phi(|x-y|))}\right)^d\, u^{-1} du \nn\\
 &\le& c_3 \frac{ a(x,y)}{|x-y|^d}
 \int^{1}_{0}
 u^{\frac{d}{2\delta_2}-1} du  \nn\\
 &\le&  c_4 \frac{ a(x,y)}{|x-y|^d}. \label{e:555}
\end{eqnarray}
Combining \eqref{e:newgfe},  \eqref{e:nup},  \eqref{e:4} and \eqref{e:555} we immediately get
part (ii) of the theorem.

\medskip
(iii)
Let $h_T (a, r)$ be defined as in \eqref{e:7.4}.
Since $a \le \Phi({\rm diam} (D)) \le (2^{-1} \wedge (2C_T)^{-2\delta_2}))T$,
we have by
\eqref{e:newgfe}--\eqref{e:4} and Lemma  \ref{l:phi_a} that
$G_D (x, y)\asymp h_T (a(x, y), |x-y|)$. The assertion then follows from
Lemma \ref{l:new1dim}.

\medskip

(iv) Note that since $d\geq 2$, we have by
\eqref{e:77} that
 \begin{eqnarray}
I &=&
\frac{\Phi(|x-y|)} {|x-y|^{d}}   \int_{\Phi(|x-y|)/T}^1  u^{-2}\left(  \frac{ |x-y|}{   \Phi^{-1} (  u^{-1} \Phi(|x-y|))} \right)^d
  \left(1\wedge \frac{  u a(x,y) }{ \Phi(|x-y|) }\right)  du\nn\\
&\leq &c_5 \frac{\Phi(|x-y|)} {|x-y|^{d}}
 \left(1\wedge \frac{   a(x,y) }{ \Phi(|x-y|) }\right) \int_0^1 u^{d/(2\delta_2)-2} \, du\nonumber \\
 &\leq &c_6 \frac{\Phi(|x-y|)} {|x-y|^{d}}
  \left(1\wedge \frac{   a(x,y) }{ \Phi(|x-y|) }\right). \label{e:5}
\end{eqnarray}
Part (iv) of the theorem  now follows from
assertion (i) of  the theorem, Lemma \ref{l:phi_a},
\eqref{e:newgfe}, \eqref{e:nup}, \eqref{e:4} and \eqref{e:5}.
\qed

\begin{cor}\label{c:gfe}
Suppose that $X$ is a one-dimensional subordinate Brownian motion  with L\'evy exponent $\Psi(\xi)=
\phi(|\xi|^2)$ with $\phi$ being a complete Bernstein function satisfying
condition {\bf (A)}.
Let $D$ be a bounded $C^{1,1}$ open subset of $\bR$
with characteristics $(R_2, \Lambda)$
and $a(x,y)=\Phi (\delta_D(x))^{1/2}
  \Phi (\delta_D(y))^{1/2}$, $x, y \in D$.

\begin{description}
\item{\rm (i)}
Suppose that for each $T>0$ there is a constant $
c_1=c_1(T, \phi)>0$ such that
\bee\label{e:6.9}
\int_r^T \frac{\Phi (s)}{s^2} \, ds \leq
c_1 \, \frac{\Phi (r)}{r}
\qquad \hbox{for }
r\in (0, T].
\eee
Then
$$
G_D(x, y)\asymp \frac{\Phi(|x-y|)} {|x-y|}
\left(1\wedge \frac{a(x,y)}{ \Phi(|x-y|)} \right).
$$

\item{\rm (ii)}
Suppose that  for every $T>0$, there is a constant $
c_2=c_2(T, \phi)>0$ so that
\bee\label{e:6.10}
\int_0^r \frac{\Phi (s)}{s^2} ds \leq
c_2 \frac{\Phi (r)}{r}
\qquad \hbox{for every } r\in (0, T].
\eee
Then
 for all  $(x, y)\in D\times D$,
$$
G_D(x, y)\asymp
  \frac{a(x,y)}{\Phi^{-1}(a(x,y))}\wedge \frac{ a(x,y)}{ |x-y|}  .
$$

\end{description}
\end{cor}

\begin{remark}\label{R:6.3}  \rm
Recall that $\delta_1,  \delta_2 \in (0, 1)$ are the constants in condition {\bf (A)}.

\begin{description}
\item{(i)} Condition \eqref{e:6.9} is satisfied when $\delta_2<1/2$.
   This is because for $t\in (0, T]$, we have
    by \eqref{e:1.9} and condition {\bf (A)},
$$
\frac{1}{\Phi (r)} \int_r^T \frac{\Phi (s)}{s^2} \, ds
= \int_r^T \frac{\phi (r^{-2})}{\phi (s^{-2})} \frac{1}{s^2} ds
\leq a_2 \int_r^T \left(\frac{s^2}{r^2}\right)^{\delta_2} \frac{1}{s^2} ds
= a_2 r^{-2\delta_2} \int_r^T s^{2\delta_2-2} ds
\leq  \frac{c}{r}.
$$

\item{(ii)}
Condition \eqref{e:6.10} is satisfied when $\delta_1>1/2$.
This is because for $t\in (0, T]$, we have
    by \eqref{e:1.9} and condition {\bf (A)},
$$
\frac{1}{\Phi (r)} \int_0^r \frac{\Phi (s)}{s^2} \, ds
= \int_0^r \frac{\phi (r^{-2})}{\phi (s^{-2})} \frac{1}{s^2} ds
\leq a_1 \int_0^r \left(\frac{s^2}{r^2}\right)^{\delta_1} \frac{1}{s^2} ds
= a_1 r^{-2\delta_2} \int_0^r s^{2\delta_1-2} ds
\leq  \frac{c}{r}.
$$
\end{description}

\end{remark}
\begin{remark}\label{R:6.4}  \rm
Let $\vp (r)={r^{1/2}}/{\phi(r)}$. Note that
\begin{eqnarray*}
a_1^{-1}\lambda^{1/2-\delta_1} \ge
\vp (\lambda r)/\vp (r)\ge a_2^{-1} \lambda^{1/2-\delta_2}\quad
&\hbox{for } \lambda \ge 1 \hbox{ and }  r \ge R_0.
\end{eqnarray*}
Let
$$\vp^* (\lambda) :=\limsup_{r \to \infty} \vp (\lambda r)/\vp (r)
\quad
\text{and}
\quad
\vp_* (\lambda) :=\liminf_{r \to \infty} \vp (\lambda r)/\vp (r).$$
Then
$$
\infty>a_1^{-1}\lambda^{1/2-\delta_1} \ge  \vp^* (\lambda) \ge    \vp_* (\lambda) \ge a_2^{-1} \lambda^{1/2-\delta_2}>0 \quad \hbox{ for } \lambda \ge 1.
$$
The upper and lower Matuszewska indices  can be computed as
$$
\alpha(\vp):=\lim_{\lambda \to \infty} (\log\vp_*(\lambda))/(\log \lambda),\quad
\beta(\vp):=\lim_{\lambda \to \infty} (\log\vp^*(\lambda))/(\log \lambda).
$$
Thus we have
$$
\infty>a_1^{-1}\lambda^{1/2-\delta_1} \ge  \vp^* (\lambda) \ge
\lambda^{\beta (\varphi)}  \ge \lambda^{\alpha (\varphi)}
\ge \vp_* (\lambda) \ge a_2^{-1} \lambda^{1/2-\delta_2}>0 \quad \hbox{ for }
\lambda \ge 1
$$
(see \cite[page 69-71]{BGT}).
Let $\wt \vp(r)=\int_{T^{-1/2}}^r {\vp (t)}/ {t}\, dt$, $ r \ge T^{-2}$.
With the change of variable $s=t^{-1/2}$, we see that
\eqref{e:6.9} is equivalent to
\bee\label{e:6.9n}
\wt \vp(r)=\int_{T^{-2}}^r\frac1{\phi (t) t^{1/2}} dt
=2\int_{r^{-1/2}}^T \frac{\Phi (s)}{s^2} \, ds \leq
c_1 \, \Phi (r^{-1/2}) r^{1/2} =\vp(r)
\qquad \hbox{for } r \ge T^{-2}.
\eee
Let $\wh \vp(r)=\int^{\infty}_r  {\vp (t)}/ {t} \, dt$.
\eqref{e:6.10} is equivalent to
\bee\label{e:6.10n}
\wh \vp(r)= \int^{\infty}_r \frac1{\phi (t) t^{1/2}} dt  =2\int_0^{r^{-1/2}} \frac{\Phi (s)}{s^2} ds \leq
c_2 \Phi (r^{-1/2}) r^{1/2}=\vp(r)
\qquad \hbox{for every } r \ge T^{-2}.
\eee
In fact, by \cite[Corollaries 2.6.2 and 2.6.4]{BGT},  $
 {\vp(r)}\asymp{\wt \vp(r)}$ for every $
  r \ge T^{-2}$
if and only if $\beta(\vp)>0$, and $
 {\vp(r)}\asymp{\wh \vp(r)}$ for every $
  r \ge T^{-2}$
 if and only if $\alpha(\vp)<0$.

By following the same proof in \cite[Section 6]{KSV5}, one can construct $\phi$ whose upper and lower
Matuszewska indices are $\alpha(\phi)=3/4$ and $\beta(\phi)=1/4$ so that
$\alpha(\vp)=1/4$ and  $\beta(\vp)=-1/4$. For such $\phi$, neither \eqref{e:6.9} nor \eqref{e:6.10} hold.
\end{remark}
\bigskip

\noindent{\bf Proof of Corollary \ref{c:gfe}.}
  Put $T=(2 \vee (2C_T)^{2\delta_2}) \Phi({\rm diam} (D))$,
 where ${\rm diam}(D)$ is the diameter of $D$ and $C_T$ is the constant in \eqref{e:77}.
For notational simplicity,
we let $a=a(x,y)$ and $r=|x-y|$.
Recall that $u_0=a/\Phi (r)$ is defined by \eqref{e:6.14}.
In view of  Theorem  \ref{t:gfe}(iii), we only need to consider $u_0 \ge1$, which
 we will assume from now on in this proof.

\medskip
\noindent (i) Since $a \le \Phi({\rm diam} (D)) \le (2C_T)^{-2\delta_2}T$,
by \eqref{e:gt2} and \eqref{e:6.9} we have
\begin{eqnarray*}
&&\frac{a}{\Phi^{-1}(a)}+\int_{r}^{\Phi^{-1}(a)} \frac{\Phi(s)}{s^2}ds
\asymp a \left(\frac{1}{\Phi^{-1}(a)}-\frac{1}{\Phi^{-1}(T)}\right)    +\int_{r}^{\Phi^{-1}(a)} \frac{\Phi(s)}{s^2}ds \\
&&=
a \int^{\Phi^{-1}(T)}_{\Phi^{-1}(a)} \frac{ds}{s^2} +\int_{r}^{\Phi^{-1}(a)} \frac{\Phi(s)}{s^2}ds \le \int_{r}^{\Phi^{-1}(T)} \frac{\Phi(s)}{s^2}ds  \le
c_1\frac{\Phi(r)} {r}.
\end{eqnarray*}
Combining this with Theorem \ref{t:gfe}(i), (iii) and Lemma \ref{l:phi_a}
establishes  part (i) of the corollary.

\medskip

\noindent (ii)
Since $a \le \Phi({\rm diam} (D))
 \le T/2$, by
\eqref{e:6.10} we have that for  $u_0 \ge1$,
\begin{eqnarray*}
 \frac{a}{\Phi^{-1}(a)}+\left(\int_{r}^{\Phi^{-1}(a)} \frac{\Phi(s)}{s^2}ds\right)^+
\asymp \frac{a}{\Phi^{-1}(a)}.
 \end{eqnarray*}
Thus by Theorem  \ref{t:gfe}(iii),
$$
G_D(x, y)\asymp  \frac{a}{\Phi^{-1}(a)}
\quad \text{ when }u_0 \ge1.
$$
Since $\Phi (r)$ is increasing in $r$, the above
 together with Theorem  \ref{t:gfe}(iii) implies that
$$
G_D(x, y) \asymp   \frac{a}{r} \wedge \frac{a}{\Phi^{-1}(a)}
 =
\begin{cases}
a/r \quad &\hbox{if } a\leq \Phi (r)  \\
a/\Phi^{-1} (a) &\hbox{if } a\geq \Phi (r).
\end{cases}
$$
\qed

We next give an example of a one-dimensional subordinate Brownian motion
with L\'evy exponent $\phi (|\xi|^2)$ that satisfies condition {\bf (A)}
but its associated function $\Phi (r)=1/\phi (r^{-2})$ satisfies neither condition
\eqref{e:6.9} nor \eqref{e:6.10}. We can get its explicit Green function estimates on bounded $C^{1,1}$ open sets by using Theorem \ref{t:gfe}(iii) but not by
Corollary \ref{c:gfe}.

 \begin{example}\label{e:g1} \rm
By the discussion at the end of the Introduction,
we know that for any $p\in \R$,
the function
$$
\phi(r)=\left(\int_0^\infty (r+t)^{-2}t^{1/2}
(\ln(t+1))^p      dt \right)^{-1}, \quad r>0
$$
is a complete Bernstein function.
Moreover, $\phi (r) \asymp r^{1/2} (\log r)^p$ when $r\geq 2$.
When $p\in [-1/2, 0)\cup (0,  1/2]$, $\phi (r)= r^{1/2} (\log (1+r))^p$ is an explicit example of a
complete Bernstein function
such that $\phi (r) \asymp r^{1/2} (\log r)^p$ when $r\geq 2$.
Let $X$ be a one-dimensional
subordinate Brownian motion with L\'evy exponent $\phi (|\xi|^2)$.
Note that
\bee \label{e:7.16}
\Phi(r)=1/ \phi(1/r^{2}) \asymp r (\log (1/r ))^{-p}
\quad \text{for }  0<r\leq 1/2 .
\eee

Let $W(x)$ be  the Lambert-W function,
that is, $W(x)$ is the unique solution to $x=W(x)e^{W(x)}$.
Suppose
$s=r (\log (1/r ))^{-p}$. Then
$$y:=p^{-1} s^{-1/p}= r^{-1/p} \log ((1/r)^{1/p})=:x e^{x}. $$
So
$\log r^{-1/p}=x=W(y)=W(p^{-1} s^{-1/p})$, which implies that
$r =\exp (-pW(p^{-1} s^{-1/p}))$. So in view of \eqref{e:1.8},
  there is a constant $c_0\in (0, 1)$ so that
\begin{equation}\label{e:LWf}
c_0 \exp (-pW(p^{-1} s^{-1/p})) \leq \Phi^{-1} (s) \leq c_0^{-1}
 \exp (-pW(p^{-1} s^{-1/p}))
\quad \hbox{ for } s \in (0, \Phi (1/2)].
\end{equation}

Recall that
$$
\int  \frac{ (\log (1/s ))^{-p}}{s} ds
=\frac1{p-1} (\log(1/s))^{1-p}+c
, \quad \text{ when } p\not=1
$$
and
$$
\int  \frac{ (\log (1/s ))^{-1}}{s} ds
=- \log(\log(1/s))+c.
$$

Suppose   $0<r \le \Phi^{-1}(a)\le 1/2$. Then
 if $p\not=1$,
\begin{eqnarray} \label{e:7.18}
 \int_{r}^{\Phi^{-1}(a)}  \frac{\Phi(s)}{s^2} ds
 & \asymp&
\int_{r}^{\Phi^{-1}(a)}  \frac{ (\log (1/s ))^{-p}}{s} ds \nn \\
&=&
\frac1{p-1}\left((\log(1/\Phi^{-1}(a)))^{1-p}-(\log(1/r))^{1-p} \right) .
\end{eqnarray}
and, if $p=1$,
\begin{eqnarray}
\int_{r}^{\Phi^{-1}(a)}  \frac{\Phi(s)}{s^2} ds
& \asymp & \int_{r}^{\Phi^{-1}(a)}  \frac{ (\log (1/s ))^{-1}}{s} ds=
\log\log(1/r)-\log\log(1/\Phi^{-1}(a))
\nn \\
 &=& \log(\log r/\log \Phi^{-1}(a)). \label{e:7.24}
\end{eqnarray}

Let $D$ be a bounded $C^{1,1}$ open set in $\R$.
We further assume that
$$
\Phi^{-1}({\rm diam} (D)) \vee {\rm diam} (D) <c_0/2,
$$
where $c_0\in (0, 1)$ is the constant in \eqref{e:LWf}.

Let $a(x,y)=\Phi (\delta_D(x))^{1/2} \Phi (\delta_D(y))^{1/2}$.
We have by Theorem  \ref{t:gfe}(iii), \eqref{e:7.16} and \eqref{e:7.18}-\eqref{e:7.24} that  for $p\not=1$,
\begin{eqnarray*}
&& G_D(x, y) \asymp
  \frac{ a(x,y)}{ |x-y|} \wedge \left( \frac{a(x,y)}
  {\Phi^{-1}(a(x,y) )}+ \left( \int_{|x-y|}^{\Phi^{-1}(a(x,y))} \frac{\Phi(s)}{s^2}ds \right)^+ \, \right) \\
&\asymp & \frac{ a(x,y)}{ |x-y|} \wedge \left( \left(
\log \frac{1}{\Phi^{-1}(a(x, y))} \right)^{-p}+
 \left( \frac{\left(\log(1/\Phi^{-1}(a (x, y)))\right)^{1-p}-
 \left(\log (1/|x-y|) \right)^{1-p}}{p-1} \right)^+ \right),
\end{eqnarray*}
 while for $p=1$,
  \begin{eqnarray*}
  G_D(x, y)
&\asymp & \frac{ a(x,y)}{ |x-y|} \wedge \left(
 \left(\log \frac{1}{\Phi^{-1}(a(x, y))}\right)^{-1}+
 \log^+ \left( \log |x-y| /\log \Phi^{-1}(a(x, y)) \right) \right),
\end{eqnarray*}
where $\log^+ x := 0 \vee \log x $. It is elementary to check that
 for $0<u,r \leq c_0/2$ and
$ c_0 u\leq v \leq c_0^{-1} u$,
 $$
 \left(\log \frac{1}{u} \right)^{-p}+
 \left( \frac{\left(\log(1/u)\right)^{1-p}-
 \left(\log (1/r) \right)^{1-p}}{p-1} \right)^+
\asymp \left(\log \frac{1}{v} \right)^{-p}+
 \left( \frac{\left(\log(1/v)\right)^{1-p}-
 \left(\log (1/r) \right)^{1-p}}{p-1} \right)^+
$$
when $p\not=1$ and
$$ \left(
\log \frac{1}{u} \right)^{-1}+
 \log^+ \left( \log r /\log u) \right)
\asymp \left( \log \frac{1}{v} \right)^{-1}+
 \log^+ \left( \log r /\log v) \right)   .
$$
Thus we have from the four displays above together with \eqref{e:LWf}
that  for $p\not=1$,
\begin{eqnarray*}
  G_D(x, y)
&\asymp & \frac{ a(x,y)}{ |x-y|} \wedge \left(
W ( p^{-1} a(x, y)^{1/p})^{-p} +
 \left( \frac{ (W ( p^{-1} a(x, y)^{1/p}))^{1-p}-
 \left(\log (1/|x-y|) \right)^{1-p}}{p-1} \right)^+ \right),
\end{eqnarray*}
 while for $p=1$,
  \begin{eqnarray*}
  G_D(x, y)
&\asymp & \frac{ a(x,y)}{ |x-y|} \wedge \Big(
 W ( a(x, y))^{-1}+
 \log^+ \left(   W (a(x, y))^{-1} \log (1/|x-y|) \right) \Big) .
\end{eqnarray*}
\end{example}

\bigskip \noindent
{\bf Acknowledgement.} The main results of
this paper, 
in particular Theorems \ref{thm:oppz}   and Corollary \ref{c:main2},
were reported
at the {\it The Sixth International Conference on Stochastic Analysis and Its Applications} held at Bedlewo, Poland, from September 9 to 14, 2012.
At the same meeting, K. Bogdan, T. Grzywny and M. Ryznar announced that
they had also obtained the factorization form estimates 
in terms of surviving probabilities and global heat kernel $p(t,c|x-y|)$
as in our Theorem \ref{thm:oppz}(i)(ii) 
for the Dirichlet heat kernels of 
a similar class
 of purely discontinuous subordinate Brownian motions considered in this paper but only for bounded $C^{1,1}$ open sets.

\vskip 0.3truein

{\bf Zhen-Qing Chen}

Department of Mathematics, University of Washington, Seattle,
WA 98195, USA

E-mail: \texttt{zqchen@uw.edu}

\bigskip

{\bf Panki Kim}

Department of Mathematics, Seoul National University,
Seoul 151-742, South Korea

E-mail: \texttt{pkim@snu.ac.kr}

\bigskip

{\bf Renming Song}

Department of Mathematics, University of Illinois, Urbana, IL 61801, USA

E-mail: \texttt{rsong@math.uiuc.edu}
\end{document}